%% file: dodec.tex
\documentclass{article} 
\usepackage{amsthm}
\usepackage{graphicx}
\usepackage{amsmath}
\usepackage{amsfonts}
\usepackage{amscd}
\usepackage{amssymb}
\usepackage{alltt}
\usepackage{makeidx}
\usepackage{txfonts}
\usepackage{url}

\input{macros}

\newtheorem{theorem}{Theorem}[section]
\newtheorem{lemma}{Lemma}[subsection]
\newtheorem{definition}[lemma]{Definition}
\newtheorem{remark}[lemma]{Remark}

\def\verbose{t}
\def\leftalign{\text{\hbox to 8in{\hss}}}

\makeindex

\begin{document}

\title{The Dodecahedral Conjecture}
\author{Thomas C. Hales\ and Sean McLaughlin \\
        \normalsize{University of Pittsburgh, Carnegie Mellon University} \\
        \normalsize{\it{hales@pitt.edu, seanmcl@gmail.com}}}
\date{}
\maketitle

\begin{abstract}  
This article gives a proof of Fejes T\'oth's Dodecahedral conjecture:
the volume of a Voronoi polyhedron in a three-dimensional packing of
balls of unit radius is at least the volume of a regular dodecahedron
of unit inradius.
\end{abstract}


\section{Introduction}

\begin{figure}[htb]
  \begin{center}
    \includegraphics[scale=0.25]{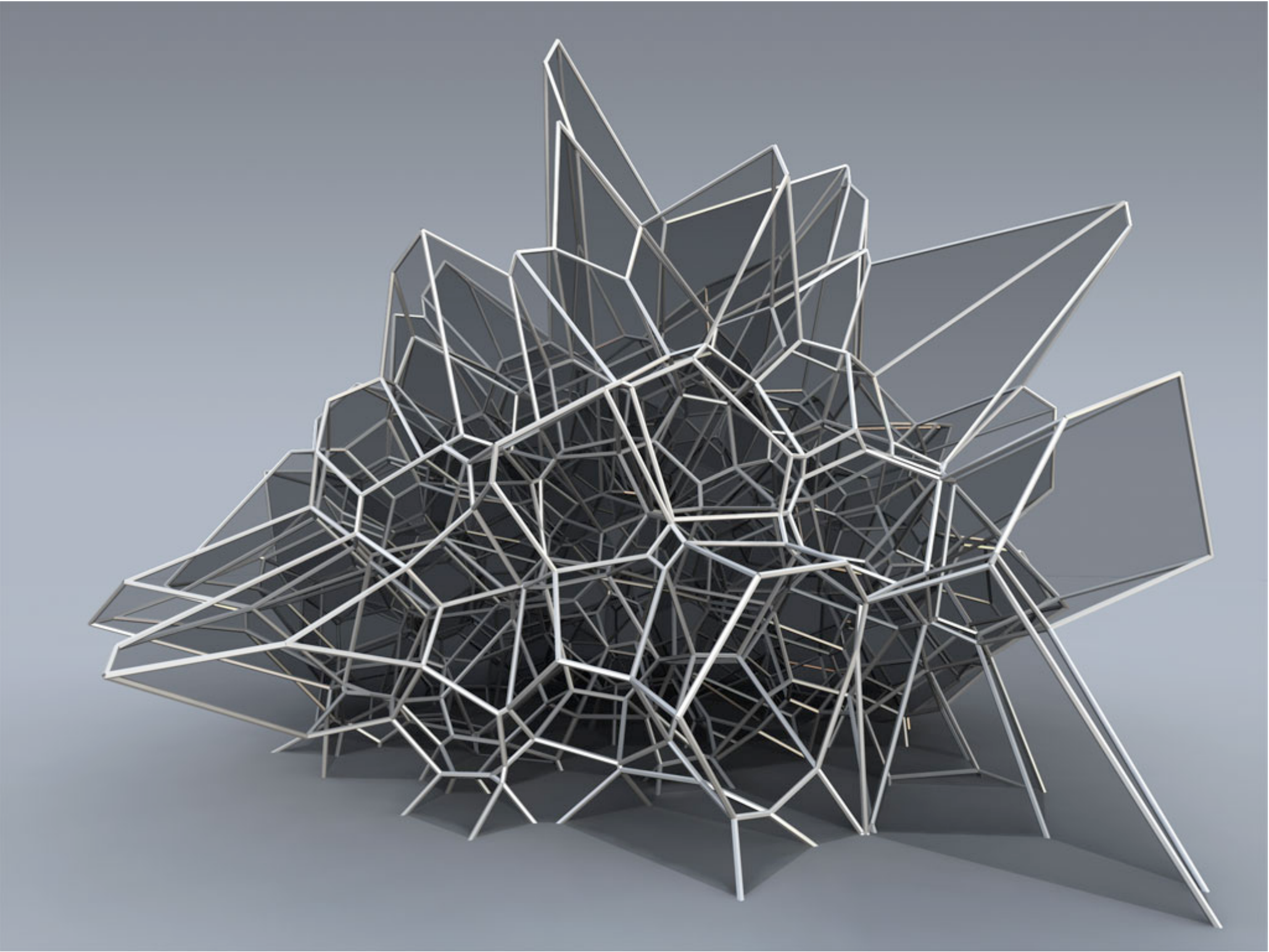}
  \end{center}
  \caption{Model of Voronoi cells by Michael Knauss and Silvan Oesterle.}
  \label{fig:voronoi}
 \end{figure}

A packing of congruent unit radius balls in three-dimensional
Euclidean space is determined by and is identified with the set
$\Lambda$ of centers of the balls.
A packing  $\Lambda$ determines a region called the Voronoi cell around
each ball.  The Voronoi cell
$\Omega(\Lambda,v)$ around a ball at $v\in \Lambda$ consists of points
of space that are closer to $v$ than to any other $w\in\Lambda$. The
Voronoi cell is a convex polyhedron containing $v$.  Figure~\ref{fig:voronoi} shows
the Voronoi cells of a finite packing.

\begin{figure}[htb]
  \begin{center}
    \includegraphics[scale=0.40]{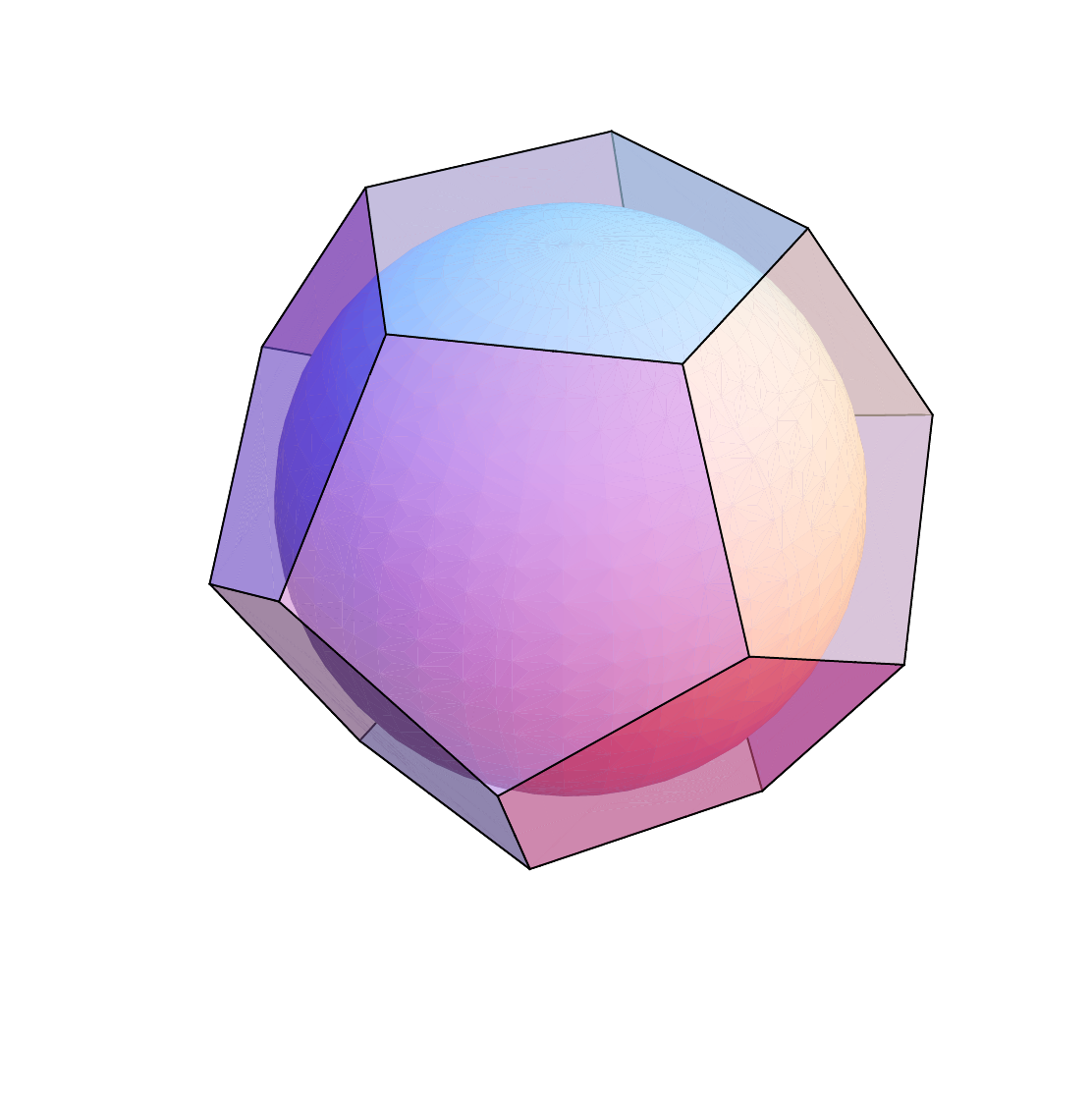}
   \end{center}
  \caption{The regular dodecahedron Voronoi cell and enclosed sphere}
  \label{fig:dod}
\end{figure}

The Dodecahedral conjecture asserts that in any packing of congruent balls of Euclidean
space every Voronoi cell has volume at least that of a regular dodecahedron
of inradius $1$ (Figure~\ref{fig:dod}).    This bound is realized by a finite 
packing $\Lambda_{dod}$
(of twelve balls and a thirteenth  at the origin) obtained
by placing a point of $\Lambda_{dod}$ at the center of each face of a regular dodecahedron (of inradius $2$).

The assertion can then be stated as the inequality
  $$
  \op{vol}(\Omega(\Lambda,v)) \ge \op{vol}(\Omega(\Lambda_{dod},0))
  $$
for every $v\in\Lambda$, and for every set of points $\Lambda\subset \ring{R}^3$
whose pairwise distances are at least the diameter $2$.
The case of equality occurs exactly when $\Omega(\Lambda,v)$ is
congruent to a regular dodecahedron of inradius $1$.

\subsection{History}

L. Fejes T\'oth made the conjecture in 1943~\cite{Toth:1943:MZ}.  
In that article, L. Fejes T\'oth sketches a proof based on an unproved
hypothesis. This hypothesis is a quantitative version of the kissing
number problem in three dimensions. This unproved hypothesis is now
generally regarded as being nearly as difficult as the Dodecahedral
conjecture itself.

L. Fejes T\'oth returned to the Dodecahedral conjecture in a number of
publications. It is a prominent part of his two books
~\cite{Toth:1972:Lagerungen}, \cite{Toth:1964:Regular}. According to
the strategy of~\cite{Toth:1972:Lagerungen}, the Dodecahedral
conjecture forms a step towards the solution of the sphere packing
problem (discussed below). In~\cite{Toth:1964:Regular} , he proved that the
Dodecahedral conjecture holds for every Voronoi cell with at most
twelve faces. This result is reviewed in Section~\ref{sec:12sphere}.
It is an ingredient in the proof presented here.

A lower bound on the volume $X$ of a Voronoi cell implies
an upper bounds on the density $4\pi/(3X)$ of packings of congruent balls
in three dimensions. The Dodecahedral conjecture gives an upper
bound on density of $0.755$.  
Upper bounds on the density based on lower bounds on the volume of a Voronoi cell
in the literature include Rogers' upper bound 0.7797~\cite{Rogers:1958:Packing}, and Muder's upper 
bounds 0.77836~\cite{Muder:1988:Voronoi} and 0.7731~\cite{Muder:1993:Bound}.

In 1993, Hsiang published what he claimed to be proofs of the Kepler
conjecture and the Dodecahedral conjecture~\cite{Hsiang:1993:IJM}.
However, the proof did not hold up to careful analysis.  ``As of this
writing, Kepler's conjecture as well as the dodecahedral conjecture
are still unproven''~\cite[p761]{Bezdek:1997:IJM}.  See also, \cite{Hales:1994:MI}.

An alternative approach to the Dodecahedral conjecture is described
in~\cite{Bezdek:1997:IJM}.  Unfortunately, a counterexample has been found to
both parts of the third conjecture of that article.  The counterexample
is described in the preprint~\cite{Hales:2002:Dodec}

K. Bezdek conjectures that the surface area of any Voronoi cell in a packing
of unit balls is at least that of a regular dodecahedron of inradius $1$.
This strengthened version of the Dodecahedral problem is still
open~\cite{Bezdek:2004:WSP}.

\subsection{The sphere packing problem}

The Kepler conjecture, also known as the sphere packing problem,
asserts that no packing of congruent balls in three dimensions has
density greater than the density of the face-centered cubic packing.
S. McLaughlin carried out the research for the proof of the
Dodecahedral conjecture at the University of Michigan
while S. Ferguson and T. Hales worked on the sphere packing problem.
Both problems were solved in 1998.

There is no strict logical connection between the two problems. The
Dodecahedral conjecture does not follow from the Kepler conjecture and
is not an intermediate step in the solution to the Kepler conjecture.
(In Fejes T\'oth's strategy, it was an intermediate step; however,
that strategy was not followed in the solution of the sphere packing
problem.) Nevertheless, the two solutions follow a similar outline and
share a significant number of methods. Both are based on long computer
calculations. Computer code was freely exchanged between the two
projects.

This article is written in a way that it is not necessary to read or
understand the solution of the sphere packing problem before reading
this article. However, for the benefit of the reader, this article
points out parallels with the sphere packing problem. It also cites
various results from that proof.

\subsection{Differences}

Although the proof of the Dodecahedral conjecture runs parallel to the
solution to the sphere packing problem, there are special difficulties
that arise in the proof of the Dodecahedral conjecture. In no sense is
it a corollary of the sphere packing problem. In the packing problem,
there turn out to be many ways to reduce the infinite ball problem to
a problem about finite clusters of balls. This multiplicity of choices
makes it possible to design many difficulties away. If one reduction
is not satisfactory one can work with another. With the Dodecahedral
problem, there is no such flexibility. The problem about finite
clusters of balls is fixed from the outset. This gives the problem a
degree of rigidity that is not present in the sphere packing problem.

\subsection{Ten years later}

Over ten years have elapsed from the completion of the research until
publication. A few words of explanation are in order. The review and
publication process for the Kepler conjecture extended from 1998 until
2006. Because of significant sharing between the Kepler conjecture and
the Dodecahedral conjecture, a ``wait and see'' attitude developed
toward the Dodecahedral conjecture. Once the Kepler conjecture was
published, the path became clear for the publication of Dodecahedral
conjecture.

The details of the proof in the current publication are essentially
the same as those in the preprint posted on the
\emph{arXiv}~\cite{Hales:1998:Dodec} in 1998. No significant errors have
emerged in the original 1998 preprint.
However,
the article has undergone various major rewrites since then, including
a significantly expanded version in 2002~\cite{Hales:2002:Dodec}.
That version is written in such a way that it is
entirely independent of the proof of the Kepler conjecture.
As an abridgment of a longer version, this article replaces 
some proofs with summaries and references.   
This article refers the reader when necessary to relevant
passages in the full version of the proof.
The computer code has
also been entirely rewritten.

A formalization project, called Flyspeck, aims to provide a complete
formalization of the proof of the Kepler
conjecture~\cite{website:FlyspeckProject, website:FlyspeckFactSheet}.
(A formal proof is one in which every logical inference of the
proof has been independently checked by computer, all the way to the
primitive axioms at the foundations of mathematics.) A parallel
project, called Flyspeck Light, aims to do the same for the proof of
the Dodecahedral conjecture. These long-term projects will take many
years to complete. Nevertheless, significant progress has already been
made toward the formal verification of the computer
code~\cite{Nipkow:2005:Tame,Obua:2008:Thesis}. The revisions in this
article incorporate the parts of Flyspeck Light that have already been
completed.

\subsection{Truncation}

The distance from the center of the regular dodecahedron to a vertex
is $t_{dod}=\sqrt{3}\tan(\pi/5)\approx 1.258$. This parameter is used
to truncate Voronoi cells; it makes volumes easier to estimate. A
similar truncation takes place in the solution to the packing problem
with truncation parameter $t_0 = 1.255$. It is a happy coincidence
that these two truncation parameters are so close to one another. A
great deal of duplicated effort might have been avoided if these two
parameters were equal. However, the parameter $t_{dod}$ cannot be
replaced with anything smaller, and although the parameter $t_0$ could
easily have been made larger, its value was already too deeply
entrenched in published articles by the time work started on the
Dodecahedral conjecture.

 As a first step towards unifying the proofs, many results can be
stated in a form that holds for all $t\in[t_0,t_{dod}]$. To transfer a
lemma from~\cite{Hales:2006:DCG} to this article, a simple process is
involved. The first step is generalization, replacing the constant
$t_0$ with a free parameter $t\in[t_0,t_{dod}]$. The second is
specialization, $t\mapsto t_{dod}$.

The number $t_0$, although rational, can be consistently treated as an
independent real transcendental in the solution to the sphere packing
problem; that is, none of the proofs involving $t_0$ rely on its exact
numerical value. The constant $t_0$ can always be replaced by a
constant $t$ in a suitably small interval about $t_0$. However, the
only way to know that this small interval is wide enough to contain
$t_{dod}$ is to study the details of the proof. We have made a
detailed study of relevant proofs in~\cite{Hales:2006:DCG} to insure
that they can be adapted to the present situation.

As a matter of terminology, a proposition for the Dodecahedral
conjecture is said to be a reparametrization of a proposition
with the constant $t_0$, if it is obtained by mechanically replacing
$t_0$ with $t_{dod}$, wherever that constant appears, and if the proof
goes through verbatim with this minor change. When this occurs, there
is nothing to further to be learned by repeating the proof.
The cited proposition already contains all the needed detail.

\subsection{Terminology}

Various notation and terminology is shared between the solution of the
sphere packing problem and this article. Vocabulary can be imported
from the sphere packing problem three different ways. The simplest way
to import a term is for the term to have precisely the same meaning in
both places. For example, the terms ``orientation'' and ``packing''
have the same meaning in both places.

The second way to import a term is by making a reparametrization of a
term that depends on the parameter $t_0$. For example, the definition
of standard component in this paper is the reparametrization of
standard component in \cite{Hales:2006:DCG}.

The third way in which terms have been imported into this proof from
the proof of the Kepler conjecture has been by structural analogy. For
example, the term {\it tame graph} is a technical notion that arises
in the solution of the sphere packing problems. At the analogous point
in this proof, a somewhat related collection of graphs appears. 
To emphasize the analogous role they play in the proof
of the Dodecahedral conjecture, they are called {\it tame Voronoi
graphs}. Similarly, there is a strong analogy between Voronoi weight
assignments in this article, and admissible weight assignments in
\cite{Hales:2006:DCG}. 

There are a some terms in that have been renamed in this version for
greater precision. {\it Special} has become {\it unstable}, {\it
standard region} has become {\it standard component}, {\it dihedral
angle} has become {\it azimuth angle}, and {\it non-external} has
become {\it internal}. {\it Planar maps} have been replaced with {\it
hypermaps}. A few terms (such as {\it distinguished}) have been
slightly redefined, when doing so is harmless.

\section{Outline}

This section gives a precise statement of the main theorem
and describes the broad outline of the proof.

\subsection{Formulation}
\label{sec:form}

This article proves the Dodecahedral conjecture in a stronger
version than that stated in the abstract. A truncation of the Voronoi
cell already has volume at least as great as that of the regular
dodecahedron. This subsection describes the truncation and states the
stronger version of the main theorem in a precise form.

Let $\Lambda$ be a packing and let $v_0\in\Lambda$. Let $B(v_0,r)$ be
a closed ball of radius $r$ centered at $v_0$. Let $\Lambda(v_0,r) =
\Lambda\cap B(v_0,r)$. Let $\CalS(\Lambda,v_0)$ be the set of all
$S=\{v_0,v_1,v_2,v_3\}\subset\Lambda(v_0,2t_{dod})$ consisting of four
distinct points such that $|v_i-v_j|\le 2t_{dod}$ for all $i,j$ and
such that the circumradius of each triangle $\{v_i,v_j,v_k\}\subset S$
is at most $\sqrt2$. Write $\op{conv}(S)$ for the convex hull of
$S\in\CalS(\Lambda,v_0)$.

Define the following truncation $\Omega_{trunc}(\Lambda,v_0)$ 
of the Voronoi cell $\Omega(\Lambda,v_0)$:
$$
\{x \in \Omega(\Lambda,v_0) \mid (x\in B(v_0,t_{dod})) \vee (\exists S\in \CalS(\Lambda,v_0).\ x \in \op{conv}(S))
$$
That is, the the Voronoi cell is truncated by intersecting it with a
ball of radius $t_{dod}$, except inside regions protected by the sets
$S\in\CalS(\Lambda,v_0)$. Note that the special packing
$\Lambda_{dod}$ satisfies 
$\Omega_{trunc}(\Lambda_{dod},0) = \Omega(\Lambda_{dod},0)$. 
The Dodecahedral conjecture takes the following strengthened form.

\begin{theorem}\label{thm:dodec}
  For every packing $\Lambda$ and every $v_0\in\Lambda$,
  $$
  \op{vol}(\Omega_{trunc}(\Lambda,v_0))\ge \op{vol}(\Omega(\Lambda_{dod},0)).
  $$
  Equality holds exactly when $\Omega(\Lambda,v_0)$ is congruent to
  $\Omega(\Lambda_{dod},0)$.
\end{theorem}

\subsection{Proof outline}
\label{sec:proof-outline}

The Lebesgue measure is translation invariant. Thus it does no harm to
assume that the center point $v_0$ of the Voronoi cell lies at the
origin: $v_0 = 0 \in \Lambda$. The assumption that $v_0=0 \in\Lambda$
remains in force for the rest of this article.

The next reduction is to replace the set $\Lambda$ with
$\Lambda(0,2t_{dod})$. This is accomplished by the following lemma,
which shows that the volumes in Theorem~\ref{thm:dodec} are
insensitive to points of $\Lambda$ outside $\Lambda(0,2t_{dod})$. The
proof appears in Lemma~\ref{lemma:trunc}.

\begin{lemma} 
  Let $\Lambda$ be any packing (with $0\in \Lambda$). Then
  $$\Omega_{trunc}(\Lambda,0) = \Omega_{trunc}(\Lambda(0,2t_{dod}),0).$$
\end{lemma}

The condition $\Lambda=\Lambda(0,2t_{dod})$ is a standing assumption for the
rest of the article.
Let $\Lambda^*=\Lambda\setminus\{0\}$ and let 
$n=\#(\Lambda^*)$ be the cardinality. The proof of
Theorem~\ref{thm:dodec} splits into two main cases: $n\le 12$ and
$n\ge 13$. In fact, L. Fejes T\'oth settles
the case $n\le 12$ 
in his book \cite{Toth:1964:Regular}. Section~\ref{sec:12sphere}
sketches Fejes T\'oth's proof.

Completely different methods treat the case when $n\ge 13$. This part
of the proof is considerably more difficult than the case treated by Fejes
T\'oth. Here is a sketch of the proof of the case $n\ge 13$. This
rough sketch will be expanded in greater detail later in
the article.

Let $\Lambda$ be a packing satisfying the standing assumptions that
$0\in\Lambda$ and $\Lambda = \Lambda(0,2t_{dod})$. The Dodecahedral
conjecture seeks a minimum to the objective function
$$
\op{vol}(\Omega_{trunc}(\Lambda,0)).
$$
This is a nonlinear optimization problem in finitely many variables.
The target value for the minimization is $\omega_{dod}=\op{vol}(\Omega(\Lambda_{dod},0))$.  When $n\ge 13$, this article proves\footnote{An examination of the proof shows the right-hand side can be improved to $\omega_{dod}+10^{-10}$.  In fact, in terms of notation to be described later in the article, the proof only
relies on the bound $\mu(\Lambda,U_F) >0$ for triangles $F$, but in fact every configuration has at least one triangle $F$ with $\mu(\Lambda,U_F) > 10^{-10}$},
$$
\op{vol}(\Omega_{trunc}(\Lambda,0))  > \omega_{dod}.
$$

Some combinatorial information about each packing $\Lambda$ is encoded
as a  graph.  
The vertex set of the graph is $\Lambda^*$.
(Because of this graph,  elements of $\Lambda^*$ are generally called vertices.)
The edge set is 
$$
E = \{\{v,w\} \mid v,w\in\Lambda^*,\quad   |v-w| \le 2t_{dod}\}.
$$
This is a planar graph.  
The Dodecahedral conjecture reduces to the case when this graph is connected.  In fact,
when the graph is not connected, this article constructs another packing
$\Lambda'$ whose graph is connected, with the same cardinality as $\Lambda$,
and such that
$$
\op{vol}(\Omega_{trunc}(\Lambda,0)) = \op{vol}(\Omega_{trunc}(\Lambda',0)). 
$$
Similarly, the conjecture reduces to the case where the graph is biconnected.
  This again involves constructing an auxiliary 
packing $\Lambda''$ of the same cardinality and whose  truncated Voronoi cell has the same volume.  Now assume that the graph of $\Lambda$ is 
biconnected.

\begin{figure}[htb]
  \begin{center}
    \includegraphics[scale=0.5]{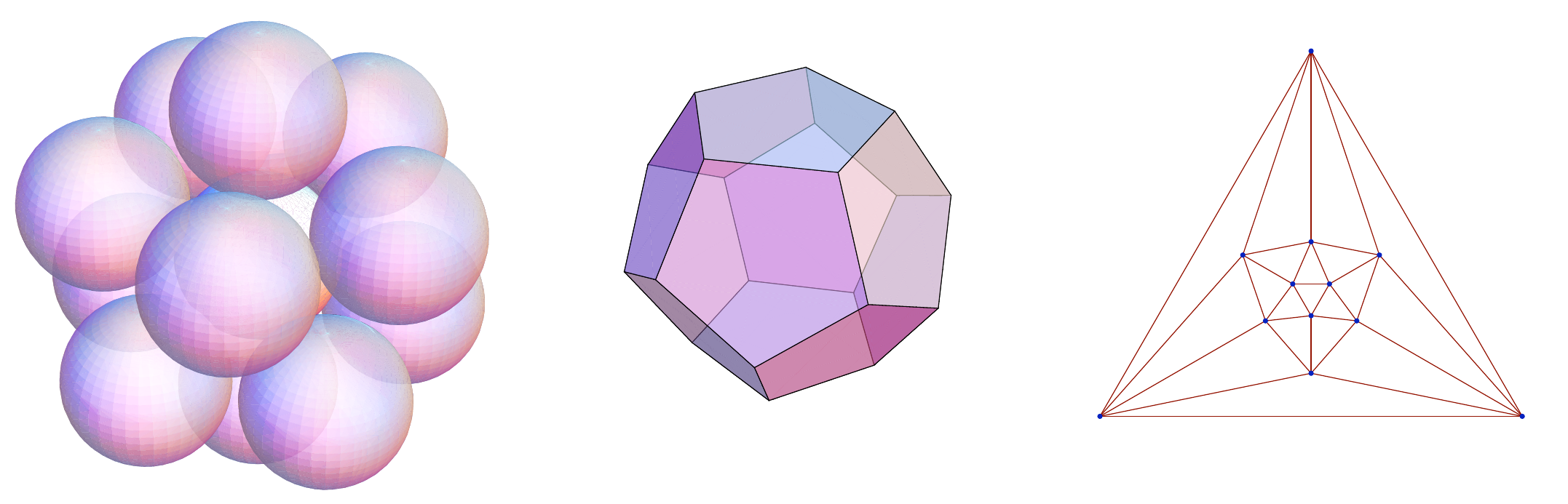}
  \end{center}
  \caption{Sphere packing $\Lambda_{dod}$, Voronoi cell $\Omega(\Lambda_{dod})$, and planar graph $G(\Lambda_{dod})$}
  \label{fig:icos}
\end{figure}

This article supposes the existence of a counterexample $\Lambda$ to
the Dodecahedral conjecture and makes a detailed study of the
properties of its graph. It defines a class of graphs (called tame
Voronoi) and proves that the graph of every counterexample is tame
Voronoi.

Tame Voronoi graphs can be described in purely combinatorial terms,
without reference to packings, Voronoi cells, and volumes. All tame
Voronoi graphs can be classified up to isomorphism. This
classification is one of the main steps of the proof. There are only
finitely many possibilities. Thus, the graph of any counterexample to
the Dodecahedral conjecture must be one of these finitely many cases.

Each tame Voronoi graph can be encoded as a hypermap $H$. (A hypermap
can be defined as a finite set, together with two permutations on that
set. The elements of the given finite set are called darts.) If $H$ is
a hypermap, let $V$ be the finite dimensional vector space of
real-valued functions on its set of darts. The pair $(H,\Phi)$ is
called a hypermap system if $\Phi$ is a set of boolean valued
functions $\phi:V^\ell \to \{\op{true},\op{false}\}$ for some $\ell$.
A hypermap system $(H,\Phi)$ is {\it feasible} if there is some
$x=(x_1,\ldots,x_\ell)\in V^\ell$ such that $\phi(x)$ holds for all
$\phi\in\Phi$. The computer code repository specifies 
a hypermap system, called
the Voronoi hypermap system, for each case $H$ arising in the
classification of tame Voronoi graphs.

Another major step of the proof is the proof that every Voronoi
hypermap system is infeasible (Theorem~\ref{thm:graph-system}). The
proof of this theorem is a case-by-case analysis based on the explicit
enumeration of tame Voronoi graphs, up to isomorphism. The feasibility
problem for each Voronoi hypermap system is converted to a system of
linear programs. The infeasibility of the Voronoi hypermap system
follows from the infeasibility of the corresponding linear program.

If there exists a counterexample to the Dodecahedral conjecture, there
is an associated Voronoi hypermap system $(H,\Phi)$. By the preceding
result, this hypermap system is infeasible. On the other hand, the
counterexample can be used to construct a feasible solution to the
system (Theorem~\ref{thm:feasible}). This contradiction shows that a
counterexample cannot exist. In this way, the Dodecahedral conjecture
is proved.

\section{Computation}

The proof of the Dodecahedral conjecture is based on a series of
computer calculations. This section briefly describes the computer
algorithms, the code implementing those algorithms, and issues of the
reliability of the computer code.

There are three main computer programs that are used in the proof. The
first is a graph generator that generates, up to isomorphism, all
planar graphs satisfying a list of properties. The second is a linear
programming package. The third is a piece of code based on interval
arithmetic that automatically proves nonlinear inequalities over the
real numbers. This section discusses each in turn.

This section also discusses some additional computer programs.
Although these computer programs, strictly speaking, are not part of
the proof, they are relevant to understanding the structure of the
proof and the reliability of the computer implementation. 
We include a brief discussion of nonlinear optimization software,
Tarski's decision procedure for real-closed fields, and formal theorem
proving packages.

\subsection{Electronic resources}

A permanent archive has been set up for all of the external resources
related to this proof~\cite{McLaughlin:2008:KeplerCode}. This archive
is under version control by Google Code~\cite{website:GoogleCode}. The
site consists of a download area where one may obtain the source code
and supporting documents to this paper. Additionally there is a
subversion~\cite{CollinsSussman:2005:Subversion} repository. This
means that the snapshot of the code and documents in the exact form
they took at the time of creation is permanently available. It
also means that any changes (for instance, a bug fix) will leave a
permanent public electronic trail. The major components
of~\cite{McLaughlin:2008:KeplerCode} are:
\begin{enumerate}
\item The source code for the three different programs used in the proof (graph generation, linear programs, and interval arithmetic inequality prover)
\item A list of all tame Voronoi graphs, up to isomorphism.
\item A list of inequalities that have been established by interval arithmetic.
\item  A list of the inequalities that have been used in linear programming.
\end{enumerate}
Additional information, such as source code documentation, is available
as well.

\subsection{Tarski arithmetic}

The Dodecahedral conjecture, after a few preliminary reductions, can
be expressed as a statement in the elementary language of the real
numbers. The elementary language of real numbers is a first order
language built from quantifiers ($\forall,\,\exists$), logical
connectives ($\land,\lor,\Rightarrow,\neg$), functions symbols for
ring operations ($+,\times,-$), variables $x_i$, and constant symbols
($0,1$). By a fundamental result of Tarski, the elementary theory of
the real numbers is decidable. Thus, the truth of the Dodecahedral
conjecture can be decided in theory by standard algorithms such as
Collin's cylindrical algebraic decomposition~\cite{Collins:1975:ATFL},
or the Cohen-H\"ormander algorithm~\cite{Hormander:1983:PDO}. However,
in practice, these decision procedures take exponential time in the
number of quantifiers, and thus are far too slow to be of practical
value for this conjecture.

To formulate the Dodecahedral conjecture as a statement in the
elementary theory of the reals, consider Voronoi cells without
truncation, centered at the origin. The volume of the regular
dodecahedron is an algebraic number $\omega_{dod}$, hence definable in
Tarski arithmetic. Also, there are a priori bounds $n$ on the number
of faces of the Voronoi cell, and thus also on the size of the
clusters of balls that give candidates for counterexamples. For
example, by adding balls to the packing to decrease the volume, the
packing becomes saturated. Assuming saturation, E. Harshbarger gives a
quick calculation of $n\le58$ faces~\cite{website:HarsbargerVoronoi}.
(The exact value of this constant is not important as long as it is
explicitly given.) The assertion of the Dodecahedral conjecture is
then expressed as an enormous conjunction of cases, with conjuncts
indexed by an explicit enumeration of all possible combinatorial
structures of a Voronoi cell, including a fixed triangulation of each
face of the cell. The fixed triangulation determines a partition of
the Voronoi cell into tetrahedra. The volume of each tetrahedron is
expressed by means of a Cayley-Menger determinant as a definable
function of the the edge lengths of the tetrahedron. Thus, for each
combinatorial structure $X$ with $m\le 58$ faces, a Tarski statement
asserts that for all vectors $v_1,\ldots,v_m\in\ring{R}^3$, if the
Voronoi cell defined by the vectors $v_i$ has combinatorial structure
$X$, then the volume of the cell is at least $\omega_{dod}$. The outer
block of $3\cdot 58$ universal quantifiers -- not to mention the
nested existential quantifiers -- is hopelessly beyond the practical
reach of current algorithms. The statement that the regular
dodecahedron is the unique minimizer can be similarly expressed.

If the entire Dodecahedral conjecture can be expressed
in Tarski arithmetic, then perhaps it is not so surprising that
many of the intermediate steps in the proof can also be so expressed.
These intermediate steps also tend to be beyond the reach of
current decision procedures.  
But here the situation is not so hopeless.  Many of these
intermediate problems involve no more than a dozen quantifiers.
One can imagine the day that these problems might fall within the
reach of decision procedures for Tarski arithmetic.

Describing various intermediate steps 
of the proof as exercises
in Tarski arithmetic is a useful point of view.  Doing so identifies
a family of subproblems that can be expressed in a common language,
and that can often be solved by similar techniques.  
The complexity of the
problems can be measured objectively by counting the number
quantifiers.

Here are some geometrical 
objects that are definable within the Tarski arithmetic
that appear in the proof of the Dodecahedral conjecture.
By Heron's formula, the circumradius $\eta(x,y,z)$ of a triangle with sides $x,y,z$ is elementary definable.
Set 
$$\eta_V(u,v,w) = \eta(|u-v|,|v-w|,|u-w|),$$ 
for three points
$u,v,w\in\ring{R}^3$.
The volume  of a tetrahedron with vertices $v_1,\ldots,v_4$ is elementary definable.  In fact, by Cayley-Menger determinants, there is a polynomial
$\Delta$ of the six squared edge lengths such that
   $$
   \sqrt{\Delta(x_{ij})}/12,\quad x_{ij}=|v_i-v_j|^2.
   $$
In particular, $\Delta(x_{ij})\ge0$, whenever the variables $x_{ij}$ have
the form $|v_i-v_j|^2$ for some points $v_i\in\ring{R}^3$~\cite[Lemma~8.1.4]{Hales:1997:DCG}.

The orientation of $v_1$ in $S=\{v_1,v_2,v_3,v_4\}$ is said to be
positive (zero, negative) 
when $v_1$ and the circumcenter of $S$ lie on the same side of (resp. on, on opposites side of)  the
plane passing through $\{v_2,v_3,v_4\}$ \cite[Lemma~5.15]{Hales:2006:DCG}.  This is an elementary condition. 

For $S=\{v_1,\ldots,v_r\}\subset\ring{R}^3$, let $\op{aff}_+(0,S)$
be the cone with apex $0$ generated by $S$:
$$
\{ t_1 v_1 + \cdots + t_r v_r \mid  t_i \ge 0\}.
$$
Write $\op{aff}_+^0(0,S)$ for the corresponding set with
strict inequality $t_i >0$.
Let $\op{conv}(S)$ be the convex hull
$$
\{ t_1 v_1 + \cdots + t_r v_r \mid t_i \ge 0,\ \sum_{i=1}^r t_i=1\}.
$$
Write $\op{conv}^0(S)$ when the inequality is strict $t_i >0$.

Define the right circular cone $\op{rcone}$ by the formula
$$\op{rcone}(v,w,h) = \{x\mid (x-v)\cdot (w-v) \ge |x-v|\,|w-v|\ h\}.$$
For a set $S$ of fixed finite cardinality, the membership conditions $x\in\op{aff}_+(0,S)$,
$x\in\op{cone}(S)$, $x\in\op{rcone}(S)$ can be expressed
in the Tarski language.

In preparing this abridged version, the proofs of numerous statements
in the Tarski language (involving a small number of quantifiers) have
been omitted.  These tend to be the arguments that can most easily
be skipped without disrupting the overall flow of the proof.
The full proofs appear at \cite{Hales:2002:Dodec}.

\subsection{Formal proof}

A formal proof is a proof in which every logical inference has been
checked, back to the foundational axioms of mathematics. Except in
trivial cases, a computer is used to generate a formal proof, because
of the large number of inferences involved. Both conventional proofs
and computer assisted proofs can be formalized. In a computer assisted
proof, this amounts to a formal verification of the correctness of the
computer code. Formal verification of computer code is a difficult
task. For that reason, formal verification tends to be reserved for
situations where correct performance is critically important, such as
the verification of aircraft control systems, cryptography algorithms,
and security protocols for the internet. There is no other means of
checking computer software that can assure reliability at levels that
remotely compare with the assurance afforded by formal verification.

There is a long term project, called Flyspeck Light, to give a formal
verification of the Dodecahedral conjecture. Although this project is
far from complete, parts of this project have already been carried
out. This means that some of the computer code for this project now
carries a proof of correctness, according to formal mathematical
standards. One such program is discussed in the next subsection.

\subsection{Graph generation}
\label{sec:gg}

The proof of the Dodecahedral conjecture is based on three separate
computer programs. The first of these is a planar graph generator. It
generates all planar graphs, up to isomorphism, that satisfy a given
list of restrictions. The restrictions include a bound on the number
of vertices in the graph, so that it is obvious that only finitely
many graphs are possible.

The correctness of the graph generator program has been the subject of
extensive mathematical investigation. Early versions of the program
were written by Hales in 1994 (in Mathematica),
1997 (in Java), and 2000 (in Java). The same computer
program is used in the Kepler conjecture and Dodecahedral conjecture.
They differ only in their input parameters. This computer program
became the subject of G. Bauer's dissertation in computer science at
the Technical University of Munich \cite{Bauer:2006:Thesis}. This
172-page dissertation translates the Java code into the formal theorem
proving system Isabelle/HOL~\cite{Paulson:1994:Isabelle} and gives a
detailed mathematical treatment of the graph theory underpinning the
computer code. The dissertation analyzes every line of code. Building
on the work of this thesis, B. Bauer and T. Nipkow have completed the
formal correctness proof of the HOL implementation of the graph
generator \cite{Nipkow:2005:Tame}. (Their published article mentions
only the Kepler conjecture, but the formal verification has been
extended to apply to the input parameters of the Dodecahedral
conjecture as well.) As a result of this work, the graph generation is
currently the most scrupulously checked part of the proof of the
Dodecahedral conjecture.

There are several published sources that provide details of this
algorithm, and there is no need to repeat details here. The basic idea
is to start with a small set of planar graphs (called `seed' graphs)
with the property that every planar graph to be classified is known to
have one of the seed graphs as a subgraph. The seed graphs are then
extended by adding one face at a time. Faces are added in all possible
ways so that it is clear at every step of the algorithm that every
biconnected planar graph will be generated. At the same time, pruning
operations discard partially completed graphs when it can be shown
that the partial graph is not a subgraph of any of the graphs to be
classified. The pruning operations prevent a combinatorial explosion
of cases. See \cite[\S5]{Hales:2003:DCG}, \cite[\S19]{Hales:2006:DCG},
\cite{Bauer:2006:Thesis}, \cite{Nipkow:2005:Tame}.

\subsection{Linear programs}

The second computer program that is used in the proof of the
Dodecahedral conjecture is linear programming. There are several
hundred linear programs that appear in the proof.

Formal verification has not yet been extended to this portion of the
computer code. However, the recent dissertation of S. Obua takes the
first steps in this direction \cite{Obua:2008:Thesis}. That work gives
a formal correctness proof of the {\it basic} linear programs that
appear in the proof of the Kepler conjecture. In particular, he has
developed all the infrastructure needed to carry out formal
correctness proofs of linear programming problems. The formalization
completed by Obua is a larger project than the formalization of the
linear programming segment of the Dodecahedral conjecture. Thus, one
can expect that the formalization of this piece of computer code will
soon follow suit.

In the 1998 proof, computer code written in C
generated the linear programs, which were then fed to the commercial
linear programming package CPLEX. In preparation for a
formal proof, the computer code has been rewritten in the programming
language Standard ML (SML), with an external
interface to the solver GLPK.

It is not necessary to trust the algorithms of the
linear programming packages (such as CPLEX and GLPK) that solve
the linear programs.  These packages produce dual certificates that
can be used to give independent verification of the solutions
of the linear programming problems.

In the proof of the Dodecahedral conjecture, the following situation
arises. The objective is to prove that the maximum of a linear
function $x\mapsto c x$ is less than a given constant $M$, when
subject to a system of linear constraints $A x \le b$. Here $x$, $c$,
$b$ are vectors with real entries and $A$ is a matrix with real
entries. The products are given by matrix multiplication of compatibly
sized matrices and vectors. A vector inequality $u\le v$ means that
$u_i \le v_i$ for each coordinate $i$. Explicit lower and upper bounds
on the variables are given: $\ell \le x \le u$. Expressed
equivalently, the objective is to show that the linear system of
inequalities
$$
A x \le b,\quad \ell \le x\le u,\quad c x \ge M
$$
has no solutions in $x$. The external linear programming package
produces a dual certificate in the form of a vector $y$, which that
package claims to have the properties
\begin{equation}
  y A = c,\quad y\ge 0,\quad y b < M.
\end{equation}
If $y$ indeed has these properties, then for any $x$ satisfying
$A x \le b$, it follows that
\begin{equation}\label{eqn:cxM}
  c x = y A x \le y b < M
\end{equation}
as desired.

Because of inexact arithmetic used by the external packages, these
identities will only be approximately correct. The imprecision in the
dual certificate can be readily eliminated as follows. If $u$ is any
vector, let $u^+$ be the vector obtained by replacing the negative
entries of $u$ with $0$, and let $u^-$ be the vector obtained by
replacing the positive entries of $u$ with $0$. By replacing the
vector $y$ with $y^+$, the vector $y$ satisfies $0\le y$. In the
following lemma, $\epsilon_1$ and $\epsilon_2$ are small error terms
that result from machine approximation. By including them in the
bounds on $c x$, a rigorous bound can be recovered.

\begin{lemma}  Suppose that the real-valued vectors and matrices
$A,A_1,A_2$, $c,c_1,c_2$, $x,b,\ell,u$, $y$ satisfy the following
relations
  $$
  A x\le b, \quad A_1 \le A \le A_2,
  \quad c_1 \le c \le c_2,\quad \ell\le x\le u,\quad
  0\le y.
  $$
Define residuals
  $$
   \epsilon_1 = c_1 - y A_2,\quad \epsilon_2 = c_2  - y A_1.
  $$
If
$$
y b + \epsilon_2^+ u^+ + \epsilon_1^+ u^- + \epsilon_2^- l^+ + \epsilon_1^- l^- < M,
$$
then $c x < M$.
\end{lemma}

\begin{proof} S. Obua has given a formal proof of this lemma in the
Isabelle/HOL system \cite[3.7.2]{Obua:2008:Thesis}. In fact, the proof
follows from a simple embellishment of Inequality~\ref{eqn:cxM}:
$$c x -y b\le c_2 x^+ + c_1 x^- -y A x= (\epsilon_2 x^+ + \epsilon_1 x^-)
  + y (A_1-A) x^+ + y (A_2 - A) x^- \le \cdots.$$
\end{proof}

The numerical data $A_1,A_2,c_1,c_2,\ell,u,y,b$ are all explicitly given,
so that the method yields explicit bounds.
It is not necessary to trust the package
that produces the certificate $y$, because
all of the assumptions of the lemma can be checked directly with
simple matrix multiplications.  The reliability of these matrix
multiplications is guaranteed by using interval arithmetic.

\subsection{Interval arithmetic}

The third major computer program that is used in the proof of
the Dodecahedral conjecture is an nonlinear-inequality prover
over the real numbers based on interval arithmetic.  
This subsection describes the
methods involved and the computer implementation.

A finite number of nonlinear functions $f_1,\ldots,f_r$ are given. It
is assumed that all functions have the same domain
$$R= [a_1,b_1] \times [a_2,b_2] \times [a_m,b_m],$$
given as a product of intervals in $\ring{R}^m$, for some $m$. The
computer program verifies that
\begin{equation}\label{eqn:fpos}
  (f_1(x) >0) \lor (f_2(x) >0) \lor \cdots \lor (f_r(x) >0),
\end{equation}
for every point $x\in R$. The approach is similar to the approach
described in R. B. Kearfott \cite{Kearfott:1996:Interval}, based on interval
arithmetic. Our methods are similar to algorithms in widespread use
for rigorous global optimization. Closely related algorithms are also
described in \cite{Zumkeller:2006:IJCAR}.

The method is based on an iteratively refined 
partition of the domain  $R$
into a finite number of smaller and smaller rectangles that cover $R$.

Start with $X=\{R\}$,  then repeat the following procedure.
Pick $T\in X$; replace $X$ with $X\setminus\{T\}$; and
calculate a lower bound $f_j(x) > a_j(T)$, 
for all $x\in T$.  If $a_j(T)\ge 0$ for some $j$, then
the desired bound (\ref{eqn:fpos}) holds on $T$.  
If the desired bound holds on $T$, then continue to the next
rectangle in $X$.
Otherwise, choose finitely many rectangles, $T_1,\ldots,T_k$ that cover $T$;
replace $X$ with $X\cup\{T_1,\ldots,T_k\}$, then repeat.
The procedure terminates when $X=\emptyset$.  (If the procedure
is applied to a false disjunction (\ref{eqn:fpos}), there is no
termination.)

When subdividing rectangles to obtain smaller covers, one does
it in such a way the the width of the rectangle tends to zero
as iteration continues.
In this way,  the lower approximations 
$a_j(T)$ to $f_j$ can be arranged to converge to the true minimum of $f_j$ on the
rectangle.

The lower bounds $a_j(T)$ to a function $f_j$ on a rectangle $T$
are obtained by methods of interval arithmetic.  The function
$f=f_j$ is generally $C^2$ on its domain.\footnote{The functions
 encountered in practice are usually $C^2$, but not always so.  When functions are not
$C^2$,  Taylor approximations are avoided.}  
The function $f$ can be expanded in a Taylor polynomial approximation with
explicit error bounds.  Derivatives are calculated by automatic
differentiation.  The error bounds are based on the Lagrange form of
the error term in the Taylor approximation.  Interval arithmetic produces rigorous bounds on the error terms.

There have been several separate implementations of the interval
arithmetic package. The source code for all of these packages is
publicly available. (The code that is used for the proof of the Kepler
conjecture is the same as the code that is used for the proof of the
Dodecahedral conjecture. Only the statements of the inequalities to be
proved differ.) The first version, written in
C++, was developed by T. Hales over the
period 1994-1998. A second version, written in C, was developed by S.
Ferguson 1995-1997. A third version, written in SML, was developed by
S. McLaughlin 2006-2008. A interval arithmetic package has also been
developed by R. Zumkeller for the theorem proving system
Coq~\cite{Bertot:2004:CoqBook}, although it has not been used to give
a formal verification of any of the inequalities that arise in the
proof of the Dodecahedral conjecture~\cite{Zumkeller:2008:Thesis}.
These implementations are all independent of one another. (Algorithms
were shared among us, but the code was independently implemented.) By
comparing the proofs of different inequalities in different systems,
we have developed a high degree of confidence that the implementations
of the algorithms are essentially correct. Of course, it would be
desirable to have a formal correctness proof, but this part of the
Flyspeck Light project has not been completed. (The SML implementation
and Zumkeller's research are partial steps in this direction.)

The list of nonlinear inequalities that are used in the proof of the
Dodecahedral conjecture appears at \cite{McLaughlin:2008:KeplerCode}.
The domains of the functions are subsets of $\ring{R}^m$, for $m\le
7$. The complexity of the verification increases rapidly with $m$.
This proof implements some of the tricks introduced in
\cite{Hales:2006:DCG} to reduce the dimension of the domain wherever
possible. Dimension reduction is based on established monotonicity
properties of the functions. (For example, the volume of a Voronoi
cell does not increase when it is intersected with a half-space.)
Whenever the functions $f_j$ are twice continuously differentiable, a
first order Taylor polynomial with explicit error bounds on the second
derivatives is used. 

The functions $f_j$ represent elementary
geometric quantities such as linear combinations of angles, dihedral
angles, solid angles, and volumes. Explicit formulas for these
functions are known involving rational functions, the square root, and
$\arctan$ functions. The typical form of a function $f_j$ is a linear
combination of terms of the form
$$
\arctan(a/\sqrt{b}),
$$
where $a,b$ are explicit polynomials on $\ring{R}^m$. In the formula
(\ref{eqn:fpos}), the number of disjuncts is usually one $r=1$, but in
some cases there are two disjuncts.
 
The computer calculations use interval arithmetic to control for
floating-point rounding errors. Every real number $x$ is represented
on the computer as an interval $[a,b]$ containing $x$, where $a$ and
$b$ are exactly representable floating point numbers
\cite{Alefeld:1983:Interval, Press:1992:NumericalRecipes}. The
calculations conform to IEEE-754 standards \cite{IEEE:1985:IEE754}.
Approximations to inverse trigonometric functions are based on
published approximations \cite{Hart:1968:Approximations}.

\subsection{Nonlinear optimization}

Previous subsections describe the three main pieces of computer code
used in the proof of the Dodecahedral conjecture: graph generation,
linear programming, and interval-arithmetic inequality proving. This
subsection describes one additional software package indirectly
involved in the proof: nonlinear-optimization. 

The disjunction of inequalities in formula (\ref{eqn:fpos}) can
be represented as a constrained minimization problem: show that
the global minimum of $f_1$ on the domain
$$
\{x\in[a_1,b_1]\times\cdots\times[a_m,b_m] \mid  f_2(x)\le 0,\ldots,
  f_r(x)\le 0 \}
$$
is positive. 

Nonlinear optimization libraries have been used to test
all the inequalities in the collection~\cite{Lawrence:1997:CFSQP},
\cite{Byrd:2006:LSNO}. 
The code generates a large random set
$X$ of points in the domain and runs the algorithm for each initial
point $x\in X$ to find a local minimum to the objective function
$f_1$. If $X$ is sufficiently large and sufficiently random, it can be
expected that one of the local minima produced to be a numerical
approximation of the global minimum.

In practice, this approach works remarkably well on this
collection of problems, largely because the functions $f_j$
tend to be rather bland from the point of view of nonlinear
optimization.  (Typically, the second derivatives of $f_j$
are small; the level surfaces of $f_j$ are approximately planar;
there are no local minima in the interior of the domain;
the global minimum occurs at a corner of the domain; and
every run of the algorithm produces the same local minimum.)
Thus, the method can usually determine the true global minimum with
high probability.

If this nonlinear optimization is not part of the proof tree, what
purpose does it serve? First of all, although we have tried to be
careful to avoid any errors in the computer code, an independent check
of the results is certainly welcome. It makes the proof more robust
against possible errors. In fact, this independent check has helped us
to spot and correct data entry errors. Secondly, the package was used
to discover inequalities that were likely to be true, and to discard
quickly inequalities that were false. The plausibly true inequalities
then became candidates for rigorous nonlinear optimization with
interval arithmetic.

\subsection{Summary}

This section has described various computer programs and algorithms
that have been used in the proof of the Dodecahedral conjecture. The
rest of the article assumes that the types of computations described
in this section can be reliably performed by computer.

Recall briefly how the three main computer programs enter into the
proof. A planar graph is associated with each potential counterexample
of the Dodecahedral conjecture. The properties of this graph is
studied, and it is shown to be a tame Voronoi graph. Using the graph
generator program, all such graphs are classified up to isomorphism.
This reduces the proof to a finite enumeration of cases. Linear
programs are then used to show that each case in this enumeration is
infeasible. The nonlinear inequalities appear in several different
parts of the proof. They are used, for example to establish that the
graph associated with a counterexample is a tame Voronoi graph.
Nonlinear inequalities are also used to justify the list of
inequalities used in the linear programs. (The linear programming
inequalities come as linear relaxations of nonlinear inequalities.)

\section{Fejes T\'oth's Reduction}\label{sec:12sphere}

L. Fejes T\'oth proved the Dodecahedral conjecture under the extra
hypothesis that $\Lambda(v_0,2t_{dod})\setminus\{v_0\}$ has at most
$12$ elements. The proof occupies about eight pages of the book
\emph{Regular Figures}~\cite{Toth:1964:Regular}. See the main theorem
of Section~33 and the main theorem of Section~41 (including the note
on page 265). Fejes T\'oth's bound is a general bound about the volume
of truncated polyhedra. The polyhedra do not need to be Voronoi cells
in a sphere packing. Here is a sketch of his proof.

\begin{theorem}
  Let $P$ be any polyhedron with at most $12$ faces that contains a
  unit sphere $S^2$. Let $B$ be the ball of radius $t_{dod}$ concentric
  with $S^2$. The volume of the intersection of $P$ with $B$ is at least
  $\omega_{dod}$. Equality holds
  exactly $P$ is congruent to the regular dodecahedron $\Omega(\Lambda_{dod},0)$.
\end{theorem}

\begin{proof}
 (Sketch) By translation, the proof reduces to the case that the
origin is the center of $S^2$. If $k<12$, there is a polyhedron with
smaller volume and $k+1$ faces, obtained by clipping a corner of the
polyhedron with a new face. Thus, assume $k=12$. View the case in
which some vertices of the polyhedron have degree greater than $3$ as
degenerate cases of polyhedra where all degrees are three, where some
of the vertices have coalesced. With these conventions, there are $12$
faces, $30$ edges, and $20$ vertices.

Let 
  $$g(x) = \begin{cases}
             \frac13 \sec^2 (x), & x \le \theta_0\\
             \frac13 \sec^2 (\theta_0), & x \ge \theta_0,
           \end{cases}$$

where $\theta_0$ is defined by $\sec(\theta_0) = t_{dod}$. For each
face $F_i$ of the polyhedron, let $w_i$ be the point
on  $S^2$ closest to the plane through
$F_i$. Let $S_i$ be the radial projection of $F_i$ to a spherical
polygon on $S^2$. The volume $\omega$ of the $t_{dod}$-truncated
polyhedron satisfies

\begin{equation}\label{eqn:V}
  \omega \ge \sum_{i=1}^{12} \int_{S_i} g(\theta(w_i,x))\,dx
\end{equation}

where $dx$ is the usual measure on $S^2$, and $\theta(w_i,x)$ is the
geodesic length of the arc on $S^2$ joining $w_i$ to $x$. The integral
on the right is exactly the volume of the truncated polyhedron
obtained by projecting each polygon $S_i$ back out to a plane through $w_i$ parallel
to $F_i$.  In particular, equality holds if the plane through each face
is tangent to $S^2$.

By the estimate of Section~33 of the book, the integral on the right is
at least
$$
120 \int_{T} g(\theta(w,x))\, dx,
$$
where $T$ is a spherical triangle with angles $\pi/2$, $\pi/5$,
$\pi/3$. Here $w$ is the vertex of the triangle $T$ that has angle
$\pi/5$. (This estimate holds more generally for any non-decreasing function
$g$.)

This integral is precisely the volume of a regular dodecahedron of
inradius $1$. Indeed, when the polyhedron is a regular dodecahedron,
the maximum of $\theta(w,x)$ over the pentagon $S_i$ is exactly
$\theta_0$. (This fact is equivalent to the definition of $t_{dod}$ as
the circumradius of the regular dodecahedron, so that no truncation
occurs.) So $g(|w-x|) = (\sec^2\theta(w,x))/3$. The inequality
(\ref{eqn:V}) is an equality. Also, each regular pentagon $S_i$ can be
triangulated into $10$ triangles with angles $\pi/5$, $\pi/2$,
$\pi/3$. The $12$ faces then give $120$ triangles $T$. The result
follows.

L. Fejes T\'oth also considers the case of equality and finds by similar
arguments that the
only minimizing polyhedron is the regular dodecahedron.
\end{proof}

\section{Geometry of Voronoi Cells}

This section describes the basic geometry of the Voronoi cell and its
truncation.

\subsection{Basic truncation}
\label{sec:in-ex}

Let $\Lambda = \Lambda(0,2t_{dod})$ be a finite packing containing
$0$, and let $\Omega(\Lambda,0)$ be the Voronoi cell. Let $B(x,r)$ be
the closed ball of radius $r$ centered at $x\in\ring{R}^3$. Let
$\Omega_0(\Lambda,0) = \Omega(\Lambda,0)\cap B(0,t_{dod})$. Also,
$\Omega_{trunc}(\Lambda,0)$ has been defined in
Section~\ref{sec:form}.

Recall $\Lambda^* = \Lambda\setminus\{0\}$.  There is a graph $G(\Lambda)$
with vertex set $\Lambda^*$, whose edges are formed by pairs
$\{v,w\}$ such that $0<|v-w|\le 2t_{dod}$.  Figure~\ref{fig:icos}
shows the regular icosahedron, which is the graph $G(\Lambda_{dod})$.

\subsubsection{One cap}

The geometry of $\Omega_0(\Lambda,0)$ will be discussed first,
then adapted to $\Omega_{trunc}(\Lambda,0)$.  The truncated cell
$\Omega_0(\Lambda,0)$ is obtained from the ball $B(0,t_{dod})$ by
removing a spherical cap

$$
C(v,t_{dod}) = \{x \in B(0,t_{dod}) \mid  |x - v| \le |x| \}
$$

\noindent for each $v\in\Lambda^*$.
Each spherical cap is bounded by a sphere of radius $t_{dod}$ and
a planar disk formed by the intersection of the bisector of 
$\{0,v\}$ with the ball $B(0,t_{dod})$.

Volume and solid angle calculations use the following functions:
\begin{equation}
\begin{array}{lll}
\phi(h,t,\lambda) &= \lambda_v t h (t+h)/6 + \lambda_s\\
A(h,t,\lambda) &= (1-h/t)(\phi(h,t,\lambda)-\phi(t,t,\lambda)),\\
\end{array}
\end{equation}
where $\lambda=(\lambda_v,\lambda_s)$.  The subscripts $v$ and $s$
on the components of $\lambda$ stand for `volume' and `solid angle.'
This terminology is justified by the following calculation.
When $\Lambda^* = \{v\}$ has a single point, the truncated Voronoi
cell is a ball of radius $t_{dod}$
with a single cap $C(v,t_{dod})$ removed.  Its volume depends
only on $|v|=2h$ (and $t_{dod}$).  
The cone at the origin generated by the points of $C(v,t_{dod})$ is
$\op{rcone}(v,h/t_{dod})$.  An elementary calculation gives the volume and
solid angle formulas.
\begin{equation}  
  \begin{array}{lll}
    \op{vol}(C(v,t_{dod})) &= \phi(h,t_{dod},(1,0))\op{sol}(\op{rcone}(v,h/t_{dod})) \\
     \phi(h,t_{dod},(0,1))&=1.\\
  \end{array}
\end{equation}
Thus, $\phi(h,t,\lambda)$ is the ratio converting  arbitrary linear combinations of volume and solid angle into solid angle.

\subsubsection{Two caps}

When $\Lambda^*=\{v,w\}$ contains two points, the two spherical
caps meet if and only if $\eta_V(0,v,w) < t_{dod}$.  (Recall that
$\eta_V$ denotes the circumradius.)   The volume formula for
$\Omega_0(\Lambda,0)$ as a function of $v,w$ is continuous
across the hypersurface $\eta_v(0,v,w)=t_{dod}$, but not analytic.
When $\eta_v(0,v,w) > t_{dod}$ the caps are disjoint and the
volume is independent of $|v-w|$, depending only on $|v|$ and $|w|$.
That is, the volume does not depend on the location of the caps,
provided they are disjoint.

When $\eta_V(0,v,w) < t_{dod}$, the volume depends on
$|v|,|w|,|v-w|$.  Note that $\eta_V(0,v,w) < t_{dod}$ implies
$|v-w|\le 2t_{dod}$, so that the graph $G(\Lambda^*)$ contains
the edge $\{v,w\}$.
By inclusion-exclusion, the volume of $\Omega_0(\Lambda,0)$ is
\begin{equation}\label{eqn:in-ex}
  \op{vol}(B(0,t_{dod})) - \op{vol}(C(v,t_{dod})) - \op{vol}(C(w,t_{dod})) +
  \op{vol}(C(v,t_{dod}) \cap C(w,t_{dod})).
\end{equation}
The set $C(v,t_{dod})\cap C(w,t_{dod})$ can be partitioned into four
regions (called quoins). 

\begin{figure}
  \begin{center}
    \includegraphics[scale=0.50]{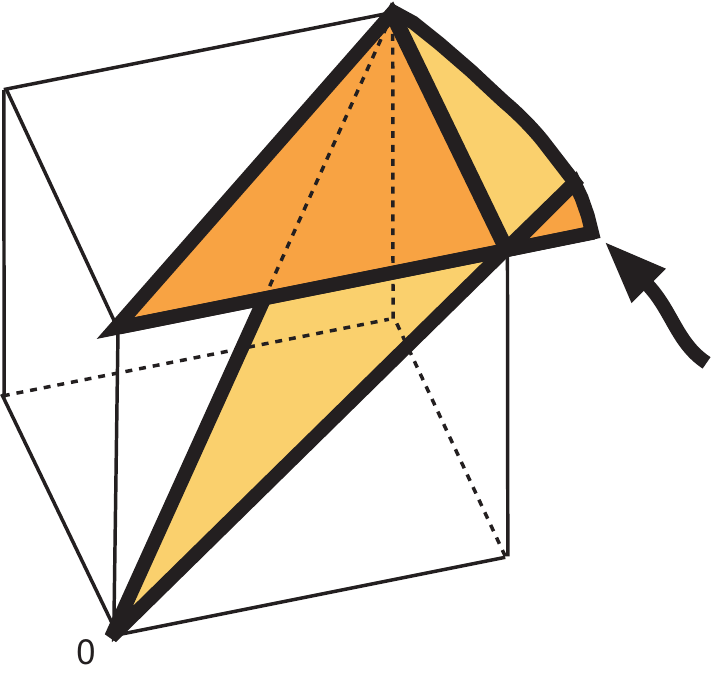}
  \end{center}
  \caption{Quoins}
  \label{fig:quoin}
 \end{figure}

Let $P_1=\op{aff}\{0,v,w\}$ be the plane through $0,v,w$. Let
$P_2$ be the plane orthogonal to $P_1$ that passes through $0$ and the
circumcenter of the triangle $\{0,v,w\}$. Then $\ring{R}^3\setminus
(P_1\cup P_2)$ contains four connected components, partitioning
$C(v,t_{dod})\cap C(w,t_{dod})$ into four quoins. Let $a=|v|/2$,
$b=\eta_V(0,v,w)$, $c=t_{dod}$. Let $q(a,b,c)$ be the volume of the
quoin given by the intersection
$$C(v,t_{dod})\cap C(w,t_{dod}) \cap H_1 \cap H_2,$$
where $H_1$ is a half-space bounded by $P_1$ and $H_2$ is
the half-space bounded by $P_2$ containing $w$.
The volume of a quoin is computed in \cite[\S7.3]{Hales:2006:DCG}.
If $a\le b\le c$, then 
\begin{equation}
\begin{array}{lll}
6\ \op{quo}(a,b,c) &= (a+2c) (c-a)^2 \arctan(e) + a(b^2-a^2) e \\
&\quad - 4c^3 \arctan(e(b-a)/(b+c)),\text{ where}\\ [0.5ex]
e^2 (b^2-a^2)&= (c^2 - b^2), \quad \text{ and } e\ge 0.
\end{array}
\end{equation}
Otherwise, $\op{quo}(a,b,c) = 0$.  By (\ref{eqn:in-ex}), this gives
the volume of $\Omega_0(\Lambda,0)$ when $\Lambda^*=\{v,w\}$.

\subsubsection{Three caps}

Lemma~\ref{lemma:4fold} shows it is not possible for four or more
spherical caps to meet.  Here we consider the case of three spherical
caps.

When $\Lambda^*=\{v_1,v_2,v_3\}$ contains three points, the three
spherical caps meet if and only if the circumradius of the simplex
$S=\{0,v_1,v_2,v_3\}$ is less than $t_{dod}$. When this happens, each
edge satisfies $|v_i-v_j|\le 2t_{dod}$, so that the graph $G(\Lambda)$ is a
triangle. Also, the circumradius of each face is at most $t_{dod} <
\sqrt2$. Thus,
$$\{0,v_1,v_2,v_3\}\in\CalS(\Lambda,0).$$
The convex hull $\op{conv}(S)$ is singled out for special truncation in
$\Omega_{trunc}$.

If the three spherical caps do not meet or meet in at most pairs,
then the volume, by inclusion-exclusion is given by a formula
similar to (\ref{eqn:in-ex}).  In particular, the volume
is given in terms of quoins, and so forth.

The case when the graph $G(\Lambda)$ is connected, but not a triangle
has particular interest.  Suppose that there is no edge between
$v_1$ and $v_3$. In this case, the volume of $\Omega_0$ depends on
$|v_i|$, $|v_1-v_2|$ and $|v_2-v_3|$, but not on $|v_1-v_3|$.
Thus, the point $v_3$ can be moved, subject to the constraints
fixing $|v_3|$ and $|v_2-v_3|$ without changing the volume.
In particular, $v_3$ can be moved until $|v_3-v_1|=2t_{dod}$.  This
results in the following simple lemma.  

\begin{lemma}\label{lemma:3tri}  
Suppose $G(\Lambda)$ contains three vertices,
is connected, but not a triangle.  Then there is another packing
$\Lambda'$ whose graph $G(\Lambda')$ is a triangle and such
that
$$
\op{vol}(\Omega_0(\Lambda,0)) = \op{vol}(\Omega_0(\Lambda',0)).
$$
\end{lemma}

The same lemma holds for $\Omega_{trunc}$ because in the situation at
hand $\Omega_0=\Omega_{trunc}$. 

\subsection{Planarity}

\begin{definition}[spherical]
For each edge $e=\{v,w\}$ of a graph $G'$ on a vertex set $\Lambda^*$, 
let $A_e$ be the
arc on the unit sphere at the origin formed by the intersection of the
sphere with $\op{aff}^0_+(0,\{v,w\})$. The graph $G'$  is said to be \emph{spherical} if the sets $A_e$ do not meet
one another, as $e$ runs over the edges of $G'$.
\end{definition}

A spherical graph is clearly a planar.

\begin{lemma}\label{lemma:planar}  
$G(\Lambda)$ is a spherical graph.
\end{lemma}

\begin{proof} 
  It suffices to show that the sets
  $\op{aff}^0_+(0,\{v,w\})$ do not meet one another, as $e=\{v,w\}$
  runs over edges.  This is a statement in Tarski arithmetic.
  A detailed proof appears in \cite[Lemma~3.2]{Hales:2002:Dodec}.  
  It is the reparametrization
  of a published theorem \cite[Lemma~3.10]{Hales:1997:DCG}.
\end{proof}

\subsection{Triangles}

This subsection describes the geometry associated with
triangles in the graph $G(\Lambda)$.  For each triangle,
there is a set $\{0,v_1,v_2,v_3\}$ such that
\begin{equation}\label{eqn:qrtet}
  |v_i| \le 2t_{dod},\quad |v_i-v_j | \le 2 t_{dod},\quad i,j\le 3.
\end{equation}

\begin{lemma} Let $S=\{0,v_1,v_2,v_3\}$ be a set of four points
that satisfies (\ref{eqn:qrtet}).  Then
there does not exist $w\in\op{conv}(S)$
that satisfies $|w-u|\ge 2$ for all
$u\in S$.
\end{lemma}

\begin{proof} This is a statement in Tarski arithmetic that can be
expressed with $12$ quantifiers ($3$ coordinates for each of
$v_1,v_2,v_3,w$). The proof is the reparametrization of
\cite[Lemma~4.15]{Hales:2006:DCG}. It is written out in full in
\cite[Lemma~3.3]{Hales:2002:Dodec}.
\end{proof}

\begin{lemma}\label{lemma:enclosed} 
Let $S=\{0,v_1,v_2,v_3\}$ be a set of four points
that satisfies (\ref{eqn:qrtet}).  Then
there does not exist $w\in\op{aff}_+(0,\{v_1,v_2,v_3\})$
(the cone with apex $0$ 
generated by positive linear combinations of $v_i$) 
that satisfies $|w-u|\ge 2$ for 
$u\in S$ and $|w|\le 2t_{dod}$.
\end{lemma}

\begin{proof} This is a statement in Tarski arithmetic. The case when
$w\in\op{conv}(S)$ is covered by the previous lemma. The remaining
case is a reparametrization of \cite[Lemma~4.19]{Hales:2006:DCG}. See
also, \cite[Cor~3.7]{Hales:2002:Dodec}.
\end{proof}

\begin{lemma}\label{lemma:4fold}  
Let $0,v_1,v_2,v_3,v_4\in\Lambda$, with
$v_i$ distinct.  The intersection $\cap_{i=1}^4 C(v_i,t_{dod})$
of spherical caps is empty.
\end{lemma}

\begin{proof}  A nonempty intersection implies that each
edge $\{v_i,v_j\}$ belongs to  $G(\Lambda)$, forming a complete graph
on four vertices.  The graph is planar by Lemma~\ref{lemma:planar}.
Its planar representation is a triangle with one vertex inside,
connected to all three vertices of the triangle.
Geometrically, this corresponds to a point $v_4\in\op{aff}_+(0,\{v_1,v_2,v_3\})$, which is impossible by the previous lemma.
\end{proof}

\begin{lemma}  Let $\{v_1,v_2,v_3\}$ and $\{v_1',v_2',v_3'\}$
be two distinct triangles in $G(\Lambda)$.  Then
$\op{aff}_+^0(0,\{v_1,v_2,v_3\})$ is disjoint from
$\op{aff}_+^0(0,\{v_1',v_2',v_3'\})$.  In particular,
$\op{conv}(\{0,v_1,v_2,v_3\})$ is disjoint from
$\op{conv}(\{0,v_1',v_2',v_3'\})$.
\end{lemma}

\begin{proof} Since the graph is planar, if the two cones
intersect, then one triangle must be contained in the other
triangle.  That is, one cone is contained in the other.  This
leads to a vertex $w$ of one triangle in the other
cone.  This is prohibited
by Lemma~\ref{lemma:enclosed}.
\end{proof}

As mentioned earlier, a three-fold intersection of spherical
caps produces a triangle in the graph $G(\Lambda)$ and a 
set 
$$\{0,v_1,v_2,v_3\}\in \CalS(\Lambda,0).$$
The next two lemmas investigate the geometry of such $\{0,v_1,v_2,v_3\}$.

\begin{lemma}\label{lemma:Q}
If $S=\{0,v_1,v_2,v_3\}$ is a set of four points such that each
face has circumradius at most $\sqrt2$ and such that
$|u-v|\le 2t_{dod}$ for $u,v\in S$, then
the circumcenter of $S$ lies in $\op{conv}(S)$.  Also, if $x
\in \op{conv}(S)$ and if some $w$ has distance at least $2$ from each
point of $S$, then $x$ is at least as close to some point of $S$ as to
$w$.
\end{lemma}

\begin{proof} This is a statement in Tarski arithmetic.
The statement about the circumcenter is \cite[Lemma~5.18]{Hales:2006:DCG}.
(No reparametrization is needed.)
If $x$  comes closer to  to $w$  than to each point of $S$, then
the Voronoi cell $\Omega(S\cup\{w\},w)$ meets $\op{conv}(S)$.
Again by \cite[Lemma~5.18]{Hales:2006:DCG}, this implies that the circumradius
of some face of $S$ is greater than $\sqrt2$, which is contrary
to hypothesis. An alternative proof 
of both parts of the lemma
is contained in \cite[Lemma~3.5,3.6]{Hales:2002:Dodec}.
\end{proof}

The following lemma justifies limiting packings to those
satisfying $\Lambda=\Lambda(0,2t_{dod})$.  Extending the packing
$\Lambda$ beyond radius $2t_{dod}$ cannot decrease the volume of the
truncated Voronoi cell.

\begin{lemma}\label{lemma:trunc}  
Let $\Lambda=\Lambda(0,2t_{dod})$.  Let $\Lambda\subset\Lambda'$,
where $\Lambda'(0,2t_{dod}) = \Lambda(0,2t_{dod})$.  Then
$$\Omega_{trunc}(\Lambda',0) = \Omega_{trunc}(\Lambda,0).
$$
That is, the truncated Voronoi cell cannot be decreased in volume by
adding additional points to the packing outside the ball $B(0,2t_{dod})$.
\end{lemma}

\begin{proof} Let $w\in \Lambda'\setminus\Lambda$ (so that  $|w|>2t_{dod}$).
Every point in $B(0,t_{dod})\cap \Omega_{trunc}(\Lambda,0)$ 
is clearly closer to $0$ than to $w$ so
belongs also to $\Omega_{trunc}(\Lambda',0)$.

Assume that $x\not\in B(0,t_{dod})$ and $x\in \Omega_{trunc}(\Lambda,0)$.
Then by the definition of $\Omega_{trunc}$,
there is some $S=\{0,v_1,v_2,v_3\}\in \CalS(\Lambda,0)$ such that
$x\in \Omega(\Lambda,0)\cap \op{conv}(S)$.  By Lemma~\ref{lemma:Q},
$x$ is closer to $0$ than to $w$.  Thus, $x\in\Omega(\Lambda',0)$.
The result follows.
\end{proof}

Assume $S=\{0,v_1,v_2,v_3\}\in \CalS(\Lambda,0)$.
From Lemma~\ref{lemma:Q}, it follows that
\begin{equation}\label{eqn:omega-3}
\op{aff}^0_+(0,\{v_1,v_2,v_3\}) \cap \Omega(\Lambda,0) = 
\op{aff}^0_+(0,\{v_1,v_2,v_3\}) \cap \Omega_{trunc}(\Lambda,0) = 
\op{conv}(S)\cap\Omega(S,0).
\end{equation}
That is, the calculation of volume can be made locally in $S$ without
reference to the position of the packing $\Lambda$.
A formula for this volume calculation appears in 
\cite[\S8.6.3]{Hales:1997:DCG}.  This gives the formula in the case
when the circumradius is at most $t_{dod}$.  When the circumradius
is at least $t_{dod}$, the spherical caps intersect in pairs,
and the inclusion-exclusion formula (using quoins)
can be used (\ref{eqn:in-ex}).
In summary, the inclusion-exclusion formula can be used for all
calculations of truncated Voronoi cells $\Omega_{trunc}(\Lambda,0)$,
except for $\{0,v_1,v_2,v_3\}\in \CalS(\Lambda,0)$ with
a circumradius less than $t_{dod}$.  In this case, 
the explicit formula just mentioned applies.

\subsection{Connecting the graph}

The volume formula for $\Omega_{trunc}(\Lambda,0)$ depends only on
$|v|$, for $v\in\Lambda^*$ and on $|v-w|$ for $\{v,w\}$ and edge of
$G(\Lambda)$. If the graph $G(\Lambda)$ is not connected, the
collections of spherical caps for two different connected components of
the graph do not intersect one another. Thus they form two or more
non-interacting ``islands'' of spherical caps that can be moved
independently around the globe $B(0,t_{dod})$ without changing the
volume of the truncated Voronoi cell. In particular, one island of
spherical caps can be moved rigidly until some vertex $v$ of on
connected component of the graph has distance exactly $2t_{dod}$ from
some vertex in another component. This connects two components of the
graph without changing volume. Thus, every truncated Voronoi cell has
the same volume of another with a connected graph. Assume without loss
of generality that the graph is connected.

A biconnected graph is defined as a connected graph that has no
articulation vertices; that is, by removing any vertex, the graph
remains connected. Consider a connected graph $G(\Lambda)$ that is not
biconnected. Then there exists some $w\in \Lambda^*$ that disconnects
the graph. Write the vertex set $\Lambda^*$ as a disjoint union of
three sets $\{w\}$, $\Lambda_1$, $\Lambda_2$ such that there is no
edge between $\Lambda_1$ and $\Lambda_2$. The volume of
$\Omega_{trunc}(\Lambda,0)$ is independent of the distances between
points of $\Lambda_1$ and $\Lambda_2$. This means that $\Lambda_2$ can
move rigidly, while constrained to preserve $|v|$ for $v\in\Lambda_2$,
$|v-u|$ for $u,v\in \Lambda_2\cup\{w\}$. Continue the rigid motion of
$\Lambda_2$ until some distance $|u-v|$ decreases to $2t_{dod}$, for
$u\in\Lambda_1$ and $v\in\Lambda_2$. (See Lemma~\ref{lemma:3tri}.)
Repeat this construction until $G(\Lambda)$ becomes biconnected. The
process does not alter the volume of $\Omega_{trunc}(\Lambda,0)$. Now
assume, without loss of generality, that $G(\Lambda)$ is biconnected.

\subsection{Standard components}

This subsection continues to assume the standing list of assumptions
on $\Lambda$. Specifically, $0\in\Lambda = \Lambda(0,t_{dod})$, and
$G(\Lambda)$ is a biconnected planar graph. Under these assumptions,
each face of $G(\Lambda)$ is a simple polygon. In particular, there
are no vertices of degree one in the graph.

For the moment, the situation can be generalized somewhat to allow
$G'$ to be any graph on the vertex set $\Lambda^*$ that is spherical and
biconnected. For each edge $\{u,v\}$ of $G'$, form the cone
$\op{aff}^0(0,\{u,v\})$. Let $X=X(G')$ be the union of these cones and
let $Y(G')$ be the complement of $X$ in $\ring{R}^3$. The open set
$Y(G')$ breaks into a finite set of connected components. Write
$[Y(G')]$ for this set of connected components. The set $[Y(G')]$ is
in natural bijection with the set of faces of the graph $G'$. Write
$U_F$ for the connected component corresponding to face $F$. The
component $U$ is said to be indexed by $F$.

Each face $F$ is identified with a sequence of vertices
$(v_1,\ldots,v_r)$ with $v_i\in G'$, giving the cyclic order of the
vertices around the face. The sequence is well-defined up to cyclic
permutation, so that $(v_1,\ldots,v_r)$ defines the same face as
$(v_r,v_1,\ldots)$. Pick the order of the cycle counterclockwise
around each face.  


Each connected component $U$ has a solid angle $\sol(U)$, which
is defined to be the area of $U\cap S^2$.  The sum of the
solid angles is the area of $S^2$:
$$
\sum_{U\in[Y(G')]} \sol(U) = 4\pi.
$$

If $U=U_F$ is a connected component and $v$ a vertex of $F$, then there
is an {\it azimuth angle} assigned to $(U,v)$ with the property
that if $F_1,\ldots,F_k$ are all the faces of $G(\Lambda)$
that contain $v$, the sum of the azimuth angles around
$v$ is $2\pi$:
$$
\sum_{i=1}^k \op{azim}(U_{F_i},v) = 2\pi.
$$
By definition, the azimuth angle equals the interior angle of the
spherical polygon $U\cap S^2$ at $v/|v|$. By Girard's formula for the
area of a triangle or polygon, for 
$F = (v_1,\ldots,v_r)$:
  \begin{equation}
  \sol(U_F) + (r-2)\pi = \sum_{i=1}^k \op{azim}(U_F,v_i)
  \end{equation}

Set $\omega(\Lambda)=\op{vol}(\Omega_{trunc}(\Lambda,0))$.
If $U$ is a connected component of $Y(G')$, set
  $$\omega(\Lambda,U) = \op{vol}(U\cap \Omega_{trunc}(\Lambda,0)).$$
Then
\begin{equation}\label{eqn:omegaU}
  \omega(\Lambda) = 
  \sum_{U\in[Y(G')]} \omega(\Lambda,U).
\end{equation}
Set $M_{dod}=0.42755$ and set  $\mu(\Lambda)= \omega(\Lambda)- 4\pi M_{dod}$.
The desired inequality can be expressed as
$\mu(\Lambda) > \mu(\Lambda_{dod})$.  
We have 
$$
  \mu(\Lambda_{dod}) \approx 0.177540.
$$
The number $\mu(\Lambda_{dod})$ is called the \emph{squander target}
in \cite{Hales:2002:Dodec}. When $U$ is a connected component of
$Y(G')$, set
\begin{equation}\label{eqn:mU}
  \mu(\Lambda,U)= \omega(\Lambda,U) - M_{dod} \sol(U).
\end{equation}
Then
$$
\mu(\Lambda) = \sum_{U\in[Y(\Lambda)]} \mu(\Lambda,U).
$$

Now specialize again to the situation where $G'=G(\Lambda)$. In this
case, write $X(\Lambda)=X(G(\Lambda))$, $Y(\Lambda)=Y(G(\Lambda))$,
and so forth. A connected component of $Y(\Lambda)$ is called a
\emph{standard component}. (This term is a reparametrization of a term
by the same name in the proof of the Kepler conjecture.)

If $U$ is indexed by a triangle $F=\{v_1,v_2,v_3\}$ in the graph,
then $U = \op{aff}_+^0(0,\{v_1,v_2,v_3\})$.  In this case,
$\omega(\Lambda,U)$ is precisely the volume of the region
already considered in Equation~\ref{eqn:omega-3}.

If $U$ is not indexed by a triangle, then
$$\omega(\Lambda,U) = \op{vol}(U\cap\Omega_0(\Lambda,0)).$$
The formula for $\omega(\Lambda,U)$ in this case follows
by inclusion-exclusion as in Section~\ref{sec:in-ex}.  Suppose
that $F$ is a face of $G(\Lambda)$ whose vertices are given
by $(v_1,\ldots,v_r)$ (listed consecutively around the face).
Set $h_i=|v_i|/2$, $t=t_{dod}$, $b^\pm_{i}=\eta_V(0,v_i,v_{i\pm 1})$,
$\lambda=(1,0)$.
By \cite[Eqn.~7.12]{Hales:2006:DCG}, inclusion-exclusion gives
\begin{equation}\label{eqn:in-ex-U}
  \begin{array}{lll}
    \omega(\Lambda,U_F) &= \sol(U_F) \phi(t,t,\lambda) + \\
    &\quad \sum_{i=1}^r (\azim(U_F,v_i) A(h_i,t,\lambda) + \quo(h_i,b_i^+,t) +
    \quo(h_i,b_i^-,t)).
  \end{array}
\end{equation}

The derivation of this formula relies on two geometric facts. First,
each quoin lies entirely in a single standard component. Second, for
each $1\le i\le r$, let $P_{\pm}$ be the open half-plane bounded by
the plane through $\{0,v_i,v_{i\pm 1}\}$, so that $P_+\cap P_- \cap V
= U_F\cap V$, for some neighborhood $V$ of $v$. Then, $C(v_i)\cap
P_+\cap P_-\subset U_F$. These facts are justified in
\cite[Lemma~12.5]{Hales:2006:DCG}. The reparametrized version
appears in \cite{Hales:2002:Dodec}.

Theorem~\ref{thm:main} shows that $\mu(\Lambda,U)$ is positive
for every standard component $U$.  The constant
$M_{dod}$ is chosen so that the minimum of $\mu(\Lambda,U)$ -- as
both $\Lambda$ and $U$  vary -- is very close to zero (about $10^{-7}$).
The function $\mu$ tends to have better numerical behavior
than $\omega$.  For that reason, even though the two
functions carry essentially the same information,
estimates are expressed in terms
of $\mu$ rather than $\omega$, whenever possible.

\begin{remark}
  \label{rem:sq} 
  The function $\mu$ is closely related to a function $\tau(\cdot,t)$
  that is used in the proof of the Kepler conjecture. The function $\mu$
  is, up to a small error term, a positive multiple of
  $\tau(\cdot,t_{dod})$. The small error term comes from the fact that
  in this article, the constant $M_{dod}$ is used, and in
  \cite{Hales:2006:DCG} the constant $M_0=1/(3 \delta_{tet})$ is used,
  where $\delta_{tet} = \sqrt8 \arctan(\sqrt2/5)$. (The constant
  $\delta_{tet}\approx 0.7796$ is Rogers's famous bound on the density
  of sphere packings.) The difference is small:
  $$M_0 - M_{dod} \approx 1.86 \times 10^{-7}.$$
  Because of the close similarity between $\mu$ and $\tau(\cdot,t_0)$,
  for every estimate involving $\tau(\cdot,t_0)$ there is apt to be an
  analogous estimate involving $\mu$. The translation involves replacing
  $M_0$ with $M_{dod}$, $t_0$ with $t_{dod}$ and rescaling the resulting
  function by an explicit positive scalar to get $\mu$.
\end{remark}

\section{The Main Estimate}

This section proves the main estimate, which gives a lower bound on
the function $\mu(\Lambda,U)$ for any standard component $U$. The
standing assumptions on $\Lambda$ remain in effect: $0\in\Lambda=
\Lambda(0,2t_{dod})$ and $G(\Lambda)$ is biconnected. This section
makes no assumptions on the cardinality of $\Lambda$ except where
explicitly stated.

\begin{theorem}\label{thm:main}  
Let $\Lambda$ be a finite packing satisfying the
standing assumptions.  Let $U_F$ be a standard component indexed by
a face $F$ of $G(\Lambda)$.  Suppose that the polygon 
$F$ has $n$ vertices.  Then
   $\mu(\Lambda,U_F) > t_n$, where 
$$
\begin{array}{lll}
 t_3 &= 0\\
 t_4 &= 0.031\\
 t_5 &= 0.076\\
 t_6 &= 0.121\\
 t_7 &= 0.166\\
 t_n &= \mu(\Lambda_{dod}),\quad n\ge 8.
\end{array}
$$
\end{theorem}

The proof of this theorem is rather long. The proof extends for
twenty pages in \cite[pp.19-38]{Hales:2002:Dodec}. The analogous
estimate in the proof of the Kepler conjecture takes a full thirty
pages \cite[pp.126-156]{Hales:2006:DCG}. We cannot pretend to give
justice to the proof under the page constraints imposed on this
version. The reader is referred to the two articles just cited for
full details of the proof. This article gives a general summary of the
ideas of the proof, with references for the reader who wishes to
pursue the proof in greater detail.

\subsection{Verifications in low dimension}

The first two cases $n=3,4$ of the theorem can be handled directly
with interval arithmetic, because they are explicit nonlinear
inequalities involving a small number of variables. The case $n=3$ can
be expressed as a nonlinear optimization problem over a tetrahedron
whose edge lengths vary in length between $2$ and $2t_{dod}$. In other
words, it is a minimization problem on the six-dimensional domain
$[2,2t_{dod}]^6$. This is readily treated by interval arithmetic
\cite{McLaughlin:2008:KeplerCode}.

The case $n=4$ can also be directly proved with interval arithmetic.
Here the optimization runs over a nine-dimensional domain. The
quadrilateral face $F=(v_1,v_2,v_3,v_4)$ is parameterized up to rigid
motion by the nine variables (3 coordinates for each of four points
minus the 3 dimensional group of rotations). Monotonicity arguments
reduce the configuration to a seven-dimensional domain. (Two of the
points $v_i$ can be rescaled $v_i \mapsto \lambda v_i$ with $0 <
\lambda \le 1$ until a constraint is met, because parallel shifts in
faces of a truncated Voronoi cell towards the origin are decreasing in
volume.) The inequality $\mu(\Lambda,U)> t_4$ on a seven-dimensional
domain can be proved directly by interval arithmetic
\cite{McLaughlin:2008:KeplerCode}.

\subsection{Strategy: superadditivity}

Define constants $D_{dod}(n,k)$ for
  $$
  n \ge 3,\quad 0\le k \le n,\quad n+k\ge 4,
  $$
by
   $$
   D_{dod}(n,k) =\begin{cases} 0.0155 & (n,k) = (3,1),\\
     t_{n+k} - D_{dod}(3,1)k & \text{otherwise}.
    \end{cases}
   $$
For $(n_1,k_1)$, $(n_2,k_2)$ in this domain,
the following superadditivity holds:
\begin{equation}\label{eqn:super}
  D_{dod}(n_1,k_1) + D_{dod}(n_2,k_2) \ge D_{dod}(n_1+n_2-2,k_1+k_2-2).
\end{equation}
In fact, this inequality follows immediately from the definitions and the
easily verified inequality,
for $m,n\ge 4$,
$$
t_m + t_n \ge t_{m+n-4} + 2 D_{dod}(3,1).
$$
Note that there are only finitely many cases involved in the
verification of this identity,
because $t_n$ is constant for $n\ge 8$.

One of the basic strategies of the proof is to give a partial triangulation
of the face $F$ (with $n$ sides)
into smaller polygons $F_1,\ldots,F_r$.   The polygon
$F_i$ will have $n_i$ sides.  Let 
$k_i$ be the number of edges $\{u,v\}$ of $F_i$ with $|u-v|\ge 2t_{dod}$.
Drawing a diagonal increases the number of oriented edges by
two, and this is the reason for the shift by two on the right
hand side of (\ref{eqn:super}).   The proof defines
a decomposition of $U_F$ into smaller components $U_i$ corresponding
to each $F_i$, gives a bound $\mu(\Lambda,U_i) > D_{dod}(n_i,k_i)$
and uses superadditivity (\ref{eqn:super}) to prove the identities:
\begin{equation}\label{eqn:mu}
  \mu(\Lambda,U_F) = \sum_{i=1}^r \mu(\Lambda,U_i).
\end{equation}
\begin{equation}\label{eqn:super-mu}
\sum_{i=1}^r \mu(\Lambda,U_i) > \sum_{i=1}^r D_{dod}(n_i,k_i)
 > D_{dod}(n,0) = t_n.
\end{equation}

The idea is that the objects $U_i$ are lower-dimensional objects
than $U_F$ (that is, the polygons  have fewer edges).
The dimension controls the complexity
of the estimates.  Thus a series of inequalities $\mu(\Lambda,U_i) > D_{dod}(n_i,k_i)$
can be expected to be easier to prove than a single inequality
$\mu(\Lambda,U_F) > t_n$ in higher dimension.  

On the other hand, if the edges of the polygons $F_i$ are allowed
to get too long, numerical experiments show that function $\mu(\Lambda,U_i)$
tends to become numerically unstable.  This prevents 
an overly aggressive  triangulation of $F$.
These experiments lead
 to a restriction of at most $3.2$ on the edge lengths.

A triple $(u,v,w)$ is called {\it unstable}
if $u,v,w$ are distinct vertices of $\Lambda^*$ such
that 
\begin{equation}\label{eqn:stable}
|u| < 2t_{dod},\ |v| < 2 t_{dod},\ |u-v|< 2t_{dod},\  |v-w|< 2t_{dod},\ \text{and } |u-w|>\sqrt8.
\end{equation}
 A pair $\{u,w\}$ is {\it
unstable} if there exists $v$ such that $(u,v,w)$ is unstable.
Otherwise it is said to be {\it stable}. Unstable edges $\{u,w\}$
create numerical instabilities and are best avoided.

\begin{remark}[strict inequalities] \label{rem:strict}
In most places in the proof, one can be sloppy about whether weak or
strict inequalities are used, but not here.  It is significant that
a weak inequality $|u - v|\ge 2t_{dod}$ is used in the definition
of the constant $k_i$.  (Deformation arguments will be used to decrease
$|u-v|$ on a closed interval and the weak inequality will keep $k_i$
constant.)
The strict inequalities  
are also significant in the definition of unstable triples.  (Deformation
arguments will be used to increase $|u|,|v|$ on the interval $[2,2t_{dod}]$
and when the upper endpoint is reached the triple becomes stable.)
\end{remark}

\subsection{Construction of subcomponents}
\label{sec:sub}

Let $F$ be a face of the graph $G(\Lambda)$ and let $U_F$ be the
corresponding standard component.  Represent $F$ as a cycle
$(v_1,\ldots,v_n)$ with $v_i\in\Lambda^*$.  The function
$\mu(\Lambda,U_F)$ depends  on $\Lambda^*$ only through
$v_1,\ldots,v_n$.  Thus, for the purpose of the proof of
Theorem~\ref{thm:main},  assume without loss of generality
that $\Lambda^* = \{v_1,\ldots,v_n\}$.

Say that $u\in\Lambda^*$ is visible from $v\in\Lambda^*$ if
$\{0,u,v\}$ is not a collinear set and if
$\op{aff}_+^0(0,\{u,v\})\subset U_F$.  When this occurs, 
call the pair $\{u,v\}$ {\it internal}.  When $\{u,v\}$
is  internal, if the edge $\{u,v\}$ is added to the
graph $G(\Lambda)$, the graph continues to be spherical.

Define $\{u,v\}$ to be a {\it distinguished} pair in $U$ if
\begin{enumerate}
\item  $|u-v|\le3.2$, 
\item $\{u,v\}$ is internal.
\item  $\{u,v\}$ is stable.
\end{enumerate}  

Inductively, build a set $X$ of distinguished edges as follows.
Start with $X=\emptyset$.
Order the distinguished pairs $\{u,v\}$ by increasing
length $|u-v|$.  Considering each
distinguished edge $\{u,v\}$ in turn,
if it satisfies
the non-crossing condition
  $$
  \op{aff}_+^0(0,\{u,v\}) \cap \op{aff}_+^0(0,\{u',v'\}) = \emptyset,
  \text{ for } \{u',v'\}\in X,
$$
then add it to $X$,

Let $G'(\Lambda)$ be the graph on vertex set $\Lambda^*=\{v_1,\ldots,v_n\}$ obtained by adding the edges $X$.  By the non-crossing
conditions, $G'(\Lambda)$ is a spherical graph.  
By the Jordan curve theorem for polygons, $Y(\Lambda)$ has two connected 
components,
$U_F$ and the complementary region $U_{F'}$.
Each
connected component $U'\subset Y(G')$ is either a subset of $U_F$ or $U_{F'}$.
Since all the added edges are internal to $U_F$, there is exactly
one component of $Y(G')$ that lies in $U_{F'}$, and that component
is equal to $U_{F'}$.  Write $[Y(G')]^* = [Y(G')]\setminus\{U_{F'}\}$
for the set of components of $Y(G')$ internal to $U_F$.

Enumerate them $U_1,\ldots,U_r$.  Define $\omega(\Lambda,U_i)$ and
$\mu(\Lambda,U_i)$, as usual,  by (\ref{eqn:mU}).  Then
(\ref{eqn:in-ex-U}) and (\ref{eqn:mu}) hold.  In fact, the justification
given for (\ref{eqn:in-ex-U}) holds verbatim in this more general
context. Let $n_i$ be the number of edges of the face
$F_i$ of $G'$ corresponding to $U_i$.  Let $k_i$ be the number of edges
$\{u,v\}$
of $F_i$ such that $|u-v|\ge 2t_{dod}$. (These are edges of $G'$ that do not
belong to $G(\Lambda)$ and edges of $G(\Lambda)$ that have
length exactly $2t_{dod}$.)  

The function $\mu(\Lambda,U_i)$ depends on $\Lambda$ only through
the vertices of $\Lambda$ on $F_i$.  Thus, for purposes of estimating
$\mu(\Lambda,U_i)$ for fixed $i$,  assume that $\Lambda$
is equal to the set of vertices of $F_i$.  The estimates can then
be expressed locally.
   This motivates the following definition.

\begin{definition}
Let $\Lambda$ be a packing such that $0\in\Lambda=\Lambda(0,2t_{dod})$.  
Let $G$ be a graph on vertex set $\Lambda^*$ 
consisting
of a single cycle containing $n\ge 3$ vertices.  
Suppose that $G$ is spherical. 
Let $U\in [Y(G)]$ be a connected component of $Y(G)$.  
The triple $(\Lambda,G,U)$ is called a {\it local configuration} if the
following conditions hold:
\begin{enumerate}
\item Every edge $\{u,v\}$ of $G$ satisfies $|u-v|\le 3.2$.
\item If
$\{u,v\}$ is internal in  $U$, then $|u-v|\ge \sqrt8$.
\item If $\{u,v\}$ is internal in $U$ and stable,
then $|u-v|\ge 3.2$.
\end{enumerate}
\end{definition}

Let $n=n(\Lambda,G,U)$ be the cardinality of $\Lambda^*$.  Equivalently,
$n$ is the number of edges in the graph $G$.  Let
$k=k(\Lambda,G,U)\le n$ be the number of edges $\{u,v\}$ of $G$ such
that $|u-v|\ge 2t_{dod}$. 

The main estimate (Theorem~\ref{thm:main}) now follows from the following
refined version of the estimate and superadditivity.

\begin{theorem}\label{thm:main'}  
Let $(\Lambda,G,U)$ be any local configuration.
Let $n=n(\Lambda,G,U)$ and $k=k(\Lambda,G,U)$ be the corresponding constants.
Assume that $n+k\ge 4$.  Then
   $$
   \mu(\Lambda,U) > D_{dod}(n,k).
   $$
\end{theorem}

\subsection{Deformations}

The proof of Theorem~\ref{thm:main'} is a total induction argument on
the cardinality $n$ of $\Lambda^*$. The induction base case is
vacuous, if the induction starts at $n=2$, since every local
configuration has $n\ge 3$. Take $n\ge 3$ and assume that
Theorem~\ref{thm:main'} holds for any local configuration with
cardinality less than $n$.

The strategy of the proof is to deform the local configuration
$(\Lambda,G,U)$ by moving a single vertex $v\in\Lambda^*$ at a time in
a way that preserves the constraint of being a local configuration,
preserves $n$, is non-increasing in $\mu(\Lambda,U)$, and is
non-decreasing in $D_{dod}(n,k)$. Under these conditions, any
counterexample to the theorem propagates to a new counterexample under
the deformation. Note that it is easily checked that $D_{dod}(n,k) \le
D_{dod}(n,k+1)$, for all $n\ge 3$ and all $n> k\ge0$. Thus, the
condition that $D_{dod}(n,k)$ is non-decreasing can be replaced with
the constraint that $k$ is non-decreasing under deformation.

The azimuth angle $\op{azim}(U,v)$ is defined for each
$v\in\Lambda^*$. Call $v$ {\it concave} in $U$ if
$\op{azim}(U,v)\ge\pi$. Otherwise, say that $v$ is convex in $U$. If
every vertex $v$ is convex in $U$, then $U$ is a convex set. When $U$
is convex, it is known that $U$ is contained in some open half-space
whose bounding plane contains the origin. Moreover, if $U$ is convex
(and $n\ge 3$), it is known that $u$ is visible from $v$ in $U$ for
any two nonadjacent vertices $u,v\in\Lambda^*$. Note that when $U$ is
convex, the conditions on local configurations require that
$|u-v|\ge\sqrt8$, for any two non-adjacent vertices $u,v\in\Lambda^*$,
because $\{u,v\}$ is automatically internal.

\subsection{Deformation at concave vertices}\label{sec:concave}

The most challenging part  of Theorem~\ref{thm:main'} is the
proof when the local configuration has a concave vertex.  This subsection
sketches the proof in that case.

The following subsections describe several different  deformations.  For each,
we describe the deformation, the starting and halting conditions on 
the deformation.  We show that the $(\Lambda,G,U)$ remains a local
configuration throughout the deformation, 
 that  $\mu(\Lambda,U)$ in non-increasing, and $k$ is
non-decreasing.

Let $v$ be concave in $U$.  Let $u,w$ be the two vertices of
$\Lambda^*$ adjacent to $v$ in $G$. 
If $|v-u|=|v-w|$, then the deformation is defined as the continuous motion
of $v$ preserving $|v|$,  moving along the bisecting plane of $\{u,w\}$
and increasing $|v-u|$.  If $|v-u|<|v-w|$, then the deformation
is defined as the continuous motion of $v$ preserving $|v|$ and
$|v-w|$ and increasing $|v-u|$.

The deformation must halt if any of the following conditions
are met.  (If the initial configuration satisfies any of these conditions,
no deformation at $v$ occurs.)
\begin{enumerate}\label{e:halt}
\item For some $v'\in\Lambda^*\setminus\{u,v,w\}$, 
$\{v,v'\}$ is internal and $|v-v'|\le \sqrt8$.
\item For some $v'\in\Lambda^*\setminus\{u,v,w\}$,
$\{v,v'\}$ is internal, stable, and $|v-v'|\le 3.2$.
\item $|v|\ge 2.2$ and $|u-v|=|v-w|=3.2$.  
\item $|v|< 2.2$, 
and $|u-v|\ge 3.07$, $|v-w|\ge 3.07$.
\end{enumerate}

By a calculation of derivatives with interval arithmetic,
the deformation is non-increasing in $\mu(\Lambda,U)$ \cite[Lemma~7.7]{Hales:2002:Dodec}.  (Although, the function $\mu(\Lambda,U)$ is potentially a 
function of a large number of variables, the derivative of $\mu$
along the deformation depends only on the six edge lengths  of
the simplex $\{0,u,v,w\}$.  This derivative calculation is within
the reach of interval methods.)
The deformation is non-decreasing in $k$, because the length
$|v-u|$ is increasing.

The deformation preserves $|v|$,
so that the constraint $0\in\Lambda=\Lambda(0,2t_{dod})$ is preserved.
For $\Lambda$ to remain a packing,  the condition
 $|u-v|\ge 2$, for $u,v\in \Lambda$ must hold.
This article does not repeat the rather technical proof that
the condition
 $|u-v|\ge 2$ is preserved, for $u,v\in\Lambda^*$.  The proof
runs a couple  of pages \cite[Lemma~7.6]{Hales:2002:Dodec}.
It is the reparametrization of \cite[Lemma~12.20]{Hales:2006:DCG}.

The next constraint is that the deformation should preserve the condition
that $G$ is spherical.  
If not, then the deformation produces
a situation where $v\in \op{aff}_+(0,\{u_1,u_2\})$ for some fixed edge $\{u_1,u_2\}$
of $G$; or $u_1\in \op{aff}_+^0(0,\{v,u\})$ for some fixed $u_1$ of $\Lambda^*$.
The first case is ruled out by \cite[Remark~p.22]{Hales:2002:Dodec}, which is the reparametrization of \cite[\S12.7,p.132]{Hales:2006:DCG}.  The second case is ruled out by 
the argument of \cite[p.27]{Hales:2002:Dodec}, a reparametrization of 
\cite[\S12.8,p.134]{Hales:2006:DCG}.

The cardinality $n$ of the vertex set $\Lambda^*$ is preserved.
The set of edges of the graph of $G$ is combinatorially determined and
remains fixed under deformation.  In particular, $G$ remains a single
cycle.  Since $G$ is spherical and consists of a single
cycle, the set $Y(G)$ has two connected components throughout the
deformation.  The component $U$ evolves continuously under deformation.   

The condition for a pair $\{v_1,v_2\}$ to be
internal is not constant under deformation.  Nevertheless,
\cite[p.23]{Hales:2002:Dodec} (or \cite[p.132]{Hales:2006:DCG}) 
shows that a pair $\{v_1,v_2\}$
cannot switch to internal when $|v_1-v_2| \le \sqrt8$.  When $\{v_1,v_2\}$ is stable,
it cannot switch to internal when $|v_1-v_2|\le 3.2$. 
The enumerated conditions on the
internal pairs $\{v_1,v_2\}$ in the definition of local configurations
now follow from the halting conditions. (In fact, the halting conditions on internal edges
can be replaced with equality, because of the
constraints on local configurations.)

The preceding arguments fully justify that the deformation preserves
the property of being a local configuration, and that
any counterexample to the lemma is propagated under
the deformation.

There is no loss in generality to assume that the first two halting conditions are never
met.  Indeed, these conditions allow a new stable internal edge $\{v,v'\}$ to be formed.  The
graph $G$ can be extended to a spherical graph $G'$ by adding the internal edge.
The component $U$ is partitioned into a disjoint union of two components $U_1,U_2$ 
of $Y(G')$ and the separating set $\op{aff}_+^0(0,\{v,v'\})$.  The vertex set
$\Lambda$ is the union of $\Lambda_1\cup\Lambda_2$ with $\Lambda_1\cap\Lambda_2=\{0,v,v'\}$
with corresponding cycles $G_1$ and $G_2$.  Both $(\Lambda_1,G_1,U_1)$ and $(\Lambda_2,G_2,U_2)$ are local configurations.  Moreover, $\mu(\Lambda,U) = \mu(\Lambda_1,U_1)+\mu(\Lambda_2,U_2)$.  By the induction hypothesis and superadditivity, the theorem follows in this case.
As the proof in this case now complete, the following arguments
assume that the first two halting conditions are never met.

Thus, the halting condition on $v$ simplifies to 
\begin{equation}\label{eqn:halt}
\begin{cases}
|u-v|=|v-w|=3.2,& \text{when } |v|\ge 2.2\\
|u-v|\ge 3.07, |v-w|\ge 3.07,& \text{otherwise.}
\end{cases}
\end{equation}
After repeating the deformation at all concave vertices, assume that the halting condition
holds for each concave vertex.  Note that the halting
condition at $v$ is incompatible with the length conditions in
the definition of an unstable triple $(v,w_1,w_2)$.  It
follows that every internal edge $\{u,w\}$
at a concave vertex $u$ is stable.  By the definition of local configuration, this
implies that $|u-w|\ge 3.2$.  Thus, the hypotheses in the following lemma are fulfilled.

\begin{lemma}\label{lemma:concave}  
Let $(\Lambda,G,U)$ be a local configuration with at least one
concave vertex.  Suppose that 
condition (\ref{eqn:halt}) holds at each concave vertex.  Suppose further  $|u-v|\ge 3.2$ whenever $\{u,v\}$ is an internal
pair with $v\in\Lambda^*$ concave and $u\in\Lambda^*$ not adjacent to $v$.
Let $n,k$ be the parameters attached
to $(\Lambda,G,U)$.  Finally, assume the induction hypothesis that Theorem~\ref{thm:main'}
holds for all $n'<n$. Then
   $$\mu(\Lambda,U) > D_{dod}(n,k).$$
\end{lemma}

\begin{proof} (Sketch)  By definitions,
$\mu(\Lambda_{dod}) = t_8 \ge t_{n+k} - D_{dod}(3,1)k = D_{dod}(n,k)$, so it is
enough to prove $\mu(\Lambda,U) > \mu(\Lambda_{dod})$. 
Define $\psi(v,\lambda)$ to be the angle opposite $\lambda$ in a triangle with sides $|v|,\lambda,t_{dod}$.  Recall the right-circular cone 
$\op{rcone}(0,v,\cos\psi(v,\lambda))$ from the discussion of Tarski arithmetic.
The proof breaks into two
cases: there are at least two concave vertices, and there is exactly one concave vertex.

Suppose that there are at least two concave vertices.  Pick two $v_1,v_2$.
Partition $U$ into three
components $U_i = U\cap \op{rcone}(0,v_i,\cos\psi(v_i,3.07/2))$, for $i=1,2$; and
$U_0= U\setminus (U_1\cup U_2)$.  It is known that
$U_1$ is disjoint from $U_2$ \cite[Lemma~3.7]{Hales:2002:Dodec}.   A  study of the geometry
of $U_1$ and $U_2$ shows the shape of $\Omega_{trunc}(\Lambda,0)\cap U_i$ and the
solid angle of $U_i$ depend only
on two parameters: $|v_i|$ and $\op{azim}(U,v_i)$.  Interval arithmetic gives
the estimates
  $$\mu(\Lambda,U_i) \ge \mu(\Lambda_{dod})/2,\quad i=1,2$$
for these two-dimensional objects \cite[\S7.2.6]{Hales:2002:Dodec}.  On the remaining piece, $\mu(\Lambda,U_0)>0$ holds
\cite[p.138]{Hales:2006:DCG}.  The sum of these terms is
  $$
  \mu(\Lambda,U) = \sum_{i=1}^3 \mu(\Lambda,U_i) > \mu(\Lambda_{dod}).
  $$

Now suppose that there is exactly one concave vertex $v$.  In this case, if $u\in\Lambda^*\setminus\{v\}$
is not adjacent to $v$, then $\{u,v\}$ is internal \cite[p.140]{Hales:2006:DCG}.  This implies
that $|u-v|\ge 3.07$ for all $u\in\Lambda^*\setminus\{v\}$.
Consider the deformation that rescales $v$ to decrease its norm $v\mapsto s v$, for 
$2/|v|\le s\le 1$. Under this deformation, $\Omega_{trunc}(\Lambda,0)\cap U$ decreases
in volume and the solid angle is unchanged, so that $\mu(\Lambda,U)$ decreases.
Combine this deformation with the deformation for concave vertices given at the beginning of this subsection so that the
constraints in the hypothesis of the lemma are preserved.  As before, the induction
hypothesis is used to avoid the first two conditions of (\ref{e:halt}).
The constants $n,k$ are unchanged; $\Lambda$ remains a packing, and so forth.
The deformation continues\footnote{A typo in \cite{Hales:2002:Dodec}  incorrectly states $|v|=2$.} until the halting condition $|v|\le 2.2$ is satisfied. 

Let $U_1 = U\cap \op{rcone}(0,v,\cos\psi(v,3.07-t_{dod}))$ and $U_0 = U
\setminus U_1$.  A study of the geometry of $\Omega_{trunc}(\Lambda,0)
\cap U_1$ shows that its volume and solid angle only depend 
on two parameters $|v|$ and
$\op{azim}(U_1,v)$.  An interval arithmetic calculation over this
two-dimensional space, using $|v|\le 2.2$, gives
$$
\mu(\Lambda,U_1) > \mu(\Lambda_{dod}).
$$
(See \cite[\S7.2.6]{Hales:2002:Dodec}. The constant $1.94159$ there  is a typo.
It should be $3.07-t_{dod}$.  The typo does not affect the proof.)
The inequality $\mu(\Lambda,U_0)>0$ holds for the same reason provided
in the case of two convex vertices.  This completes
the proof of Lemma~\ref{lemma:concave}.

\end{proof}

\subsection{Deformation at convex vertices}

The results of the previous subsection
reduce the proof of Theorem~\ref{thm:main'} to
the case where $U$ is convex at every vertex.  A convex
spherical polygon on a unit sphere has perimeter
at most $2\pi$.  A polygon $F$ in $G(\Lambda)$ projects to
a spherical polygon on the unit sphere $S^2$.  An edge $\{u,v\}$ of
$F$ satisfies bounds $|u|,|v|\in[2,2t_{dod}]$, $|u-v|\ge 2$. This
implies that every edge of the spherical polygon has arc length
at least $\theta=2\arcsin(1/(2t_{dod}))$, and that the number
of edges is at most seven ($2\pi/\theta < 8$).

In this subsection, the geometry is much more explicit than in the
previous subsection, because $U$ is convex and $F$ has at most seven sides.
Let $(G,\Lambda,U)$ be a local configuration.
Let $(n,k)$ be the associated parameters.
This section describes the proof of Theorem~\ref{thm:main'} under 
the total induction hypothesis on $n$ and assuming the truth
for local configurations already treated.
The method, again, is to produce a deformation of the local
configuration by moving one vertex $v$ at a time.

Let $(\Lambda,G,U)$ be a local configuration with $U$ convex.
Let $v\in\Lambda^*$.  
Let $u,w$ be the two vertices of $\Lambda^*$
adjacent to $v$.  By convexity, every pair $\{u,v'\}$,
with $v'\in\Lambda^*\setminus\{u,v,w\}$, is internal.

\subsubsection{First convex deformation}

Consider the deformation that fixes $|v|$ and $|v-w|$ and
moves $v$ to decrease $|v-u|$.  
The deformation halts (or never starts) once any
of the following conditions holds.
\begin{enumerate}\label{e:halt-convex}
\item For some $v'\in\Lambda^*\setminus\{u,v,w\}$, 
$|v-v'|\le \sqrt8$.
\item For some $v'\in\Lambda^*\setminus\{u,v,w\}$,
the pair $\{v,v'\}$ is  stable and $|v-v'|\le 3.2$.
\item $\op{azim}(U,v)\ge\pi$.
\item There exists $v'\ne v$ such that 
$(u,v',w)$ is an unstable triple and $|u-w|\le3.2$.
\item $|v-u|=\sqrt8$.
\item $|v-u|=2t_{dod}$.
\item $|v-u|=2$.
\end{enumerate}

As with deformations at nonconvex vertices, 
the deformation of a local
configuration remains a local configuration.  Here, in the
convex situation, the proof is more elementary, because the
geometry is explicit.  For instance, the condition that
$\Lambda$ remains a packing follows immediately from the halting
conditions, because every nonadjacent vertex gives an internal
edge.
The function $\mu(\Lambda,U)$ is non-increasing under the
deformation by \cite[Lemma~7.8]{Hales:2002:Dodec}.  The value of
$k$ in non-decreasing by the halting conditions. 

\subsubsection{Second convex deformation}

Let $(\Lambda,G,U)$ be a local configuration with $U$ convex.
Consider the deformation that fixes $|u-v|$ and $|u-w|$ and
moves $v$ to increase or decrease $|v|$.  The direction
of the deformation is chosen to decrease $\mu(\Lambda,U)$.
The deformation halts (or never starts) once any
of the following conditions holds.
\begin{enumerate}\label{e:halt-convex2}
\item For some $v'\in\Lambda^*\setminus\{u,v,w\}$, 
 $|v-v'|\le \sqrt8$.
\item For some $v'\in\Lambda^*\setminus\{u,v,w\}$,
the pair $\{v,v'\}$ is  stable, and $|v-v'|\le 3.2$.
\item $\op{azim}(U,v)\ge\pi$.
\item $|u-v|\not\in\{ 2,2t_{dod},\sqrt8\}$.
\item $|u-w|\not\in \{2,2t_{dod},\sqrt8\}$.
\item $|v|=2$.
\item $|v|=2t_0$.
\end{enumerate}

By an interval arithmetic calculation of derivatives,
the function $\mu(\Lambda,U)$ does not have a local
minimum, provided none of the halting conditions hold
\cite[Lemma~7.10]{Hales:2002:Dodec}.  That is, the deformation can always continue
to decrease $\mu(\Lambda,U)$ until a halting condition is met.

\subsubsection{Completion of the proof}

\begin{proof} (Sketch) With these two deformations at hand, the proof
of Theorem~\ref{thm:main'} can be completed. Let $n$ be the
cardinality of $\Lambda^*$. If $n=3$, the inequality of the theorem is
an inequality in six variables and can be verified directly by
interval arithmetic \cite{McLaughlin:2008:KeplerCode},
\cite[\S7.4.1]{Hales:2002:Dodec}. Now assume that $n>3$. Furthermore,
by previous estimates, $n\le 7$.

By induction, it may be assumed that
the theorem is established for all $n'<n$.  By previous arguments,
it may be assumed that the theorem is known for $U$ with a concave
$v$ (and the same value of $n$).

By the induction argument and the reduction to the convex case,
there is no loss in generality to assume that the first three
halting conditions (for both deformations) never occur.  

The  halting condition (4) of the first convex deformation is rather strange:
There exists $v'\ne v$ such that 
$(u,v',w)$ is an unstable triple and $|u-w|<3.2$.
(It was needed in the interval arithmetic verifications
that prove the monotonicity of $\mu(\Lambda,U)$.)  Note
that $v$ and $v'$ are both adjacent to $u,w$ when this
halting condition holds.  This implies that $n=4$.  

Consider the case $n=4$.  The dimension of a general configuration
is nine, parameterized by four lengths $|v|$ for $v\in\Lambda^*$,
four lengths $|u-v|$ for edges $\{u,v\}$ of $G(\Lambda)$, and $|u-v|$
for one internal pair $\{u,v\}$.
Even
if the halting condition (4) becomes binding at $v$, deformations
can continue at the three other vertices, until some halting
condition holds at each vertex.  Eventually the deformations
reduce the dimension of the configuration to at most three.
Interval arithmetic finishes this case~\cite[\S7.4.2]{Hales:2002:Dodec}.

With the case $n=4$ out of the way, the
the halting condition for the first convex deformation
reduces to
  $$
  |u-v| \in\{ 2, 2t_{dod}, \text{ or } \sqrt8\}.
  $$
The first convex deformation can be applied at each vertex 
so that this condition holds for every edge $\{u,v\}$ of $G(\Lambda)$.
Then, the halting conditions (4) and (5) of the second convex
deformation now never occur.  The second convex deformation can
be applied until $|v|=2$ or $|v|=2t_{dod}$ for each $v\in\Lambda^*$.

The local configuration $(\Lambda,G,U)$ is now a low-dimensional
object.  The only remaining continuous parameters are the lengths of
$(n-3)$ internal pairs $\{u,v\}$ needed to triangulate the
$n$-gon $G(\Lambda)$.  Thus, it has dimension $m=n-3$,
for $5\le n\le 7$.

Unfortunately, the proof does not end here with a simple interval
arithmetic calculation in low dimension.  It does not end
here because there is no control on the
lengths of the triangulating diagonals, and without any such
control the calculations are simply too numerically unstable.

A different strategy completes the proof, based on truncated
corner cells.  Although the dimension of this problem is now small,
this argument requires several pages.
See \cite[pp.30-38]{Hales:2002:Dodec}.  It is modeled on a published
$14$-page argument in the solution to the sphere packing problem
\cite[\S\S13.2-13.11]{Hales:2006:DCG}, following the same strategy.
The final parts of this section give a brief summary of the two principal  methods that are used.

\subsubsection{Dealing with unstable edges}
If $(u,v,w)$ is an unstable triple, then the
halting conditions force $|u|=|w|=|u-v|=|v-w|=2$.  (This
relies on the strictness of the inequalities (\ref{eqn:stable})
defining stability,
as mentioned in Remark~\ref{rem:strict}.)
The value of $|v|$ can be $2$ or $2t_{dod}$. In this final
stage of the proof, contrary to the constraints of Section~\ref{sec:sub}
on distinguished pairs,
it is now permitted to split the region $U$ into two pieces
$U_F,U_{F'}$ separated by $\op{aff}_+^0(0,\{u,w\})$, where
the pair $\{u,w\}$ is unstable and $|u-w|\le 3.2$.  
One of these pieces is indexed
by a triangle $F$.  Because of instability, the usual inequality
$\mu(\Lambda,U_F)> D_{dod}(3,1)$ does not hold.  To compensate,
stronger inequalities are proved for $\mu(\Lambda,U_{F'})$.
The deformations can continue on the component $U_{F'}$. Nevertheless,
it must be remembered that the induction hypothesis and the 
reduction to the convex case do not cover the stronger 
inequality for $\mu(\Lambda,U_{F'})$ that is now needed.

\subsubsection {Truncated corner cells}  If $U$ has no unstable
internal $\{u,w\}$, then the following argument gives the
desired bound.  The component $U$ is partitioned into $n+1$
parts, one $U_v$ for each vertex $v\in\Lambda^*$ and a final part
for the remainder $U_0 = U\setminus(\cup_{v\in\Lambda^*} U_v)$.
The function $\mu(\Lambda,U)$ is a sum of terms $\mu(\Lambda,U_v)$
and $\mu(\Lambda,U_0) >0$.  The function $\mu(\Lambda,U_v)$
is a function of the six edges of $\{0,v,u,w\}$ (with $u,w$
adjacent to $v$) and most of these edges are fixed in length
by the deformations.  So the function $\mu(\Lambda,U_v)$
is readily bounded with interval arithmetic.  Each part
$$\Omega_{trunc}(\Lambda,0)\cap U_v$$
is called a {\it truncated corner cell}.  The defining conditions
for $U_v$ are
  $$
  U_v = U \cap \op{rcone}(0,v,\cos\psi(v,1.6)) \cap H(0,v,u) \cap H(0,v,w),
  $$
where $\psi(v,\lambda)$ is the angle defined in Section~\ref{sec:concave},
and $H(0,v,v')$ is the open half-space containing $v$, bounded by
the plane through $0$ and through the circumcenter
of $\{0,v,v'\}$, orthogonal to the plane of $\{0,v,v'\}$.  The
sets $U_v$, as $v$ ranges over $\Lambda^*$, are disjoint from one
another.
We refer the reader to the unabridged version of the proof for details.
\end{proof}

\section{Classification of Tame Hypermaps}

This section turns to the problem of classifying a large finite collection
of planar graphs. For combinatorial simplicity, this classification is phrased
in terms of hypermaps, which are defined in the first subsection.
The next subsection shows how a sphere packing $\Lambda$ gives a hypermap.
The rest of the section is devoted to the classification problem.
A final subsection shows how a counterexample $\Lambda$ to the Dodecahedral
conjecture gives one of the hypermaps classified in this section.

\subsection{Hypermap}

A hypermap is a tuple $(D,e,n,f)$, where $D$ is a finite
set, and $e,n,f$ are three permutations on that set that
compose to the identity:
$e\circ n\circ f = I$.  The elements of $D$ are called darts.
The permutations $e,n,f$ are called the edge permutation,
node permutation, and face permutation, respectively.
(A hypermap was previously defined as a finite set $D$ with
two permutations $f,n$, which amounts to the same thing,
since $e$ is uniquely determined by $f,n$.)

If $m$ is any permutation on $D$, write $D/m$ for the
set of orbits in $D$ under $m$.  Similarly, if $G$ is any
group of permutations on $D$, write $D/G$ for the set
of orbits of $D$ under $G$.  In particular, $D/\tangle{e,n,f}$
is the set of orbits under the group generated by $e,n,f$.
An orbit of $D$ under $f$ ($n$, or $e$) is called a face (resp.
node, or edge).

\begin{figure}[htb]
  \begin{center}
    \includegraphics[scale=0.50]{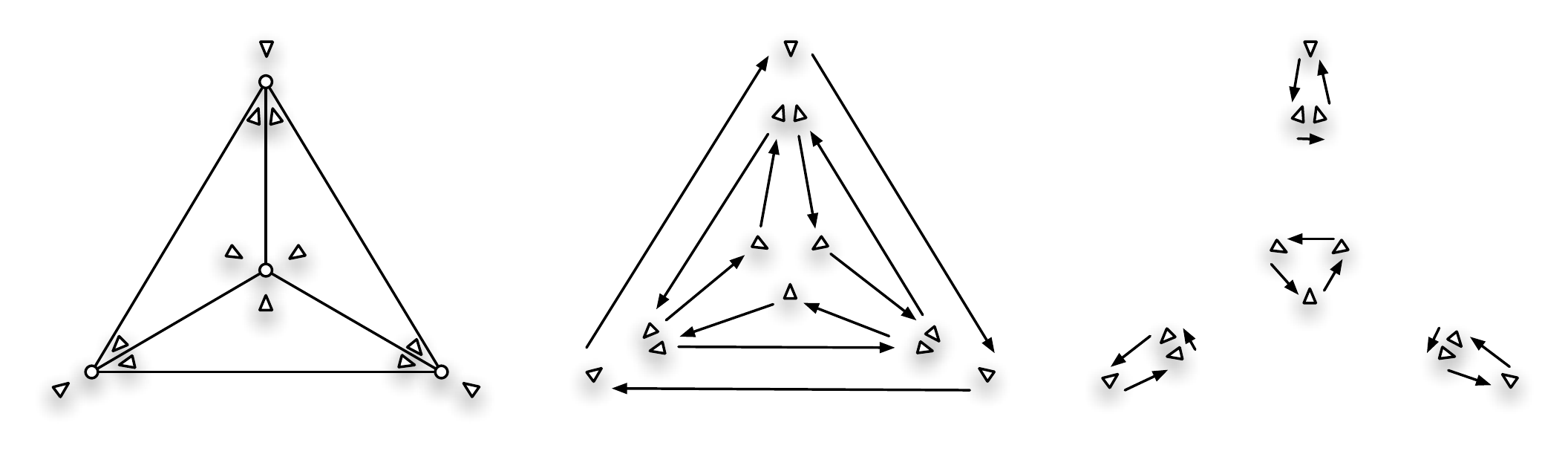}
   \end{center}
  \caption{A planar graph as hypermap, with faces and nodes}
  \label{fig:hypermap}
\end{figure}

A planar graph gives a hypermap
by the following procedure.  Starting with a planar graph,
place a dart at each angle (or equivalently at the tail of each
directed edge).  That is, at a vertex of degree $k$,
place $k$ darts, one between each consecutive pair of edges.
The face permutation has a cycle for
each face of the planar graph and  traverses the
darts in a counterclockwise direction around each face.
The node permutation has a cycle for each vertex and traverses
the darts in a counterclockwise direction around each vertex.
The edge permutation is defined by the relation $e\circ n\circ f=I$.
It can be interpreted as an involution that pairs a dart associated
with the tail of a directed edge with a dart associated with the tail
of the oppositely directed edge..
See Figure~\ref{fig:hypermap}.

Hypermaps are the primary combinatorial object used by Gonthier in the
formalization of the Four-Color theorem in
Coq~\cite{Gonthier:2005:FourColor}. Hypermaps, by being purely
combinatorial, are more convenient to represent on a computer than
planar graphs.

Not all hypermaps arise from a planar graph in this way.
Those that do have two special properties.  They are involutive
and planar in the following sense.   The definition
of planar hypermap is the standard condition on the Euler
characteristic, translated into the language of hypermaps.

\begin{definition}\label{def:involutive}
\mbox{}
\begin{itemize}
\item The hypermap $(D,e,n,f)$ is involutive, if $e$ is an involution:
$$
 e^2 = I.
$$
\item The hypermap $(D,e,n,f)$ is planar, if
   $$
   \#(D/e) + \#(D/n) + \#(D/f) = \# D + 2\#(D/\tangle{e,n,f}).
   $$
where $\# X$ denotes the cardinality of $X$.
\end{itemize}
\end{definition}

\subsection{Packings and hypermaps}\label{sec:ph}

Let $\Lambda$ be a packing satisfying $0\in\Lambda=\Lambda(0,2t_{dod})$.  The graph
$G(\Lambda)$ is planar.   Assume that $G(\Lambda)$ is biconnected. 
This subsection describes in greater detail the
hypermap $H(\Lambda)$ attached to $G(\Lambda)$.

For each $v\in \Lambda^*$, let 
  $$E(v) = \{u \mid \{v,u\} \text{ is an edge in the graph } G(\Lambda)\}.$$
For each $u\in E(v)$, there is a half-plane $P_+(v,u)$ containing $u$, bounded by
the line through $\{0,v\}$.  There is a cyclic order on the half-planes $P_+(v,u)$,
moving in a counterclockwise circle around the ray emanating from $0$ through $v$.
Write $\sigma_v$ for the cyclic permutation on $E(v)$, given by this ordering.

Define the set of darts by
$$
D = \{(0,v,u,\sigma_v u) \mid u\in E(v)\}.
$$
Define face, edge, and node permutations on $D$ by
$$
\begin{array}{lll}
  f(0,v,w,u) &= (0,w,\sigma_w^{-1} v,v),\\
  e(0,v,w,u) &= (0,w,v,\sigma_w v),\\
  n(0,v,w,u) &= (0,v,u,\sigma_v u).
\end{array}
$$
A formal calculation shows that $(D,e,n,f)$ is an involutive hypermap.

The nodes of $(D,e,n,f)$ are in bijection with $\Lambda^*$ under the correspondence:
$$
  (0,v,w,u) \mapsto v\in\Lambda^*.
$$
The edges of $(D,e,n,f)$ are in bijection with the edges of $G(\Lambda)$ under the
correspondence:
$$
  (0,v,w,u) \mapsto \{v,w\}.
$$
The faces of $(D,e,n,f)$ are in bijection with the faces of $G(\Lambda)$ under the
correspondence:
$$
  x = (0,v,w,u) \mapsto F=((f^{-1}x)_3,x_3,(f x)_3,(f^2 x)_3,\ldots),
$$
where $y_3$ is the third component of the four-tuple $y$ and each face of $G(\Lambda)$ is represented
as usual as cycle $(v_1,\ldots,v_n)$.
The set of darts is in bijection with the set of oriented edges of $G(\Lambda)$ under
the correspondence:
$$
(0,v,w,u)\mapsto (v,w).
$$
The graph $G(\Lambda)$ is connected.  This implies that $\tangle{e,n,f}$ acts
transitively on $D$.

The graph $G(\Lambda)$ is planar.  If $V,E,F$ are the number of vertices, edges,
and faces, then $V-E+F=2$; or equivalently, $V+E+F = 2 + 2E$.
Under the bijections just described, this implies that the cardinalities
of these sets satisfy
$$
     \#(D/e) + \#(D/n) + \#(D/f) = \# D + 2\#(D/\tangle{e,n,f}).
$$
Thus, the hypermap $(D,e,n,f)$ is planar.

\begin{remark}
In \cite{Hales:2002:Dodec}, the basic combinatorial structure is called a {\it planar map} rather
than hypermap.  In that article, the combinatorial structure is represented in
computer code as a finite set of faces
 $$
  \{F_1,F_2,\ldots,F_r\},
 $$
and each face is represented as a cycle $(v_1,\ldots,v_n)$ of vertices.  This is
essentially equivalent to a hypermap.  This representation
is converted to a hypermap by sending $(v_1,v_2,\ldots,v_n)$ to the dart
$(0,v_1,v_2,\sigma_{v_1} v_2)$.  There are $n$-choices of which vertex $v_i$ to list
first in the cycle, and by taking all choices, $n$ darts are obtained.  Running
through all faces in this way, all darts are constructed.  In the opposite direction,
an earlier argument describes how a face of the hypermap gives a face of $G(\Lambda)$,
expressed as a cycle.
\end{remark}

\subsection{Tameness}

This subsection defines a collection of hypermaps called tame Voronoi
hypermaps. The classification of these hypermaps, up to isomorphism,
is one of the main steps of the proof of the Dodecahedral conjecture.
Bauer and Nipkow have formally proved the classification in
Isabelle~\cite{Nipkow:2005:Tame}.  This Isabelle proof was originally
designed for the classification of the sphere packing problem. Nipkow
and McLaughlin have modified that proof to cover the Dodecahedral
conjecture as well. The modified Isabelle proof is found
at~\cite{McLaughlin:2008:KeplerCode}.

Let $H=(D,e,n,f)$ be a hypermap.  A face of $H$ (that is, an orbit
of $D$ under $f$) is said to be a triangle, quadrilateral, pentagon, etc. if the cardinality
of the orbit is $3$, $4$, $5$, respectively.  Two nodes are said
to be adjacent
if there is an edge $\{x,y\}$ of $H$ such that $x$ belongs to one of the nodes and
$y$ belongs to the other.

  Let $v$ be a node of $H$ (that is, an orbit
of $D$ under $n$).  
A node $v$ is said to have type $(p,q,r)$, if the cardinality of $v$ is $p+q+r$ and if
there are $p$ triangular faces, $q$ quadrilateral faces, and $r$ other faces
that share a dart with $v$.  The cardinality of a node is also
called its degree.

\begin{table}[h]
  \begin{center}
    \begin{tabular}{|c|c|c|c|c|c|c|} 
      \hline
      $b(p,q)$ & 0 & 1 & 2 & 3 & 4 \\
      \hline
      0 & * & * & * & 0.093 & 0.125  \\
      1 & * & * & 0.092 & 0.093 & *  \\
      2 & * & 0.133 & 0.062 & * & *  \\
      3 & * & 0.043 & 0.118 & * & *  \\
      4 & 0.053 & 0.051 & * & * & *  \\
      5 & 0.004 & * & * & * & *  \\
      6 & 0.121 & * & * & * & * \\
      7 & * & * & * & * & * \\
      \hline
    \end{tabular}
  \end{center}
  \caption{Vertex types}
  \label{vertexTable}
\end{table}

Define constants $b(p,q)$ by Table~\ref{vertexTable}.
If $(p,q)$ falls outside this table, or if the entry is marked $*$, then
set $b(p,q)=\mu(\Lambda_{dod})$.
Let $t_n$, $n\ge3$, be the collection of constants defined in Theorem~\ref{thm:main}.

A weight assignment of a hypermap $(D,e,n,f)$ is a function $w:D \to \ring{R}$ that
is constant on faces: $w(f x) = w(x)$ for $x\in D$.  A weight assignment is said to be a Voronoi weight assignment
if the following properties hold:
\begin{enumerate}
\item If the face containing $x$ has cardinality $m$, then $w(x)> t_m$.
(In particular, $w(x)>0$ for all $x$.)
\item Let $F\subset D$ 
be any face with cardinality $m \ge 5$.  Let $y\in F$.
Let $V$ be a set of nodes, each meeting $F$, such
that no two are adjacent to one another.  
Assume that the type of each node of $V$ is $(4,0,1)$.
Let $X =(\cup V)\setminus F$;
that is, the set of darts in nodes in $V$ except those in $F$.
Let $m'$ be the cardinality of $V$.
Then
$$
w(y) + \sum_{x\in X} w(x) > t_m  +  0.016\, m'.
$$
\item If the node of $x$ has type $(p,q,0)$ and degree $m=p+q$, then
  $$
  \sum_{i=1}^m w(n^i x) > b(p,q).
  $$
\end{enumerate}
The total weight of a weight assignment $w$ is defined to be
$$
\sum_{x\in [D/f]}  w(x),
$$
where $[D/f]$ is a set of representatives of the orbits of $D$ under $f$.

\begin{figure}[htb]
  \begin{center}
    \includegraphics[scale=0.50]{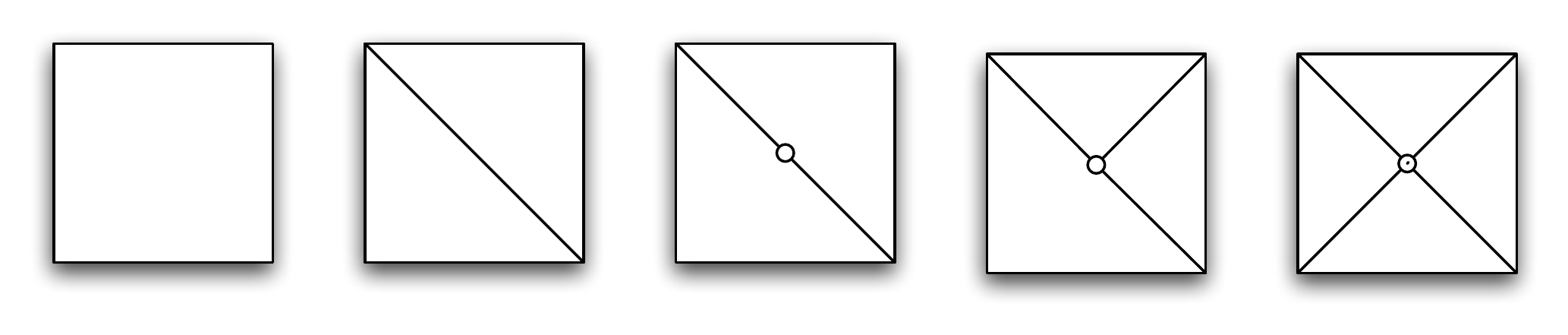}
   \end{center}
  \caption{Quad types}
  \label{fig:quadtype}
\end{figure}

\begin{definition}
$H=(D,e,n,f)$ is said to be a tame Voronoi hypermap if the following conditions
hold.
\begin{enumerate}

\item $H$ is an involutive, planar hypermap.
\item $H$ is connected.  
That is, $D$ is a single orbit under $\tangle{e,n,f}$.
\item (Simple face) Every face of $H$ meets every node of $H$ in at most
dart.

\label{def:tame6}
\item The number of nodes is at least $13$.

\item The cardinality of each face of $H$ is at least $3$ and at most $7$.
\label{def:tame1}

\item 
\label{def:tame2}
(Triangle types) 
Let $v_1,v_2,v_3$ be any three nodes of $H$ such that $v_i$ is adjacent to $v_j$
for each $i\ne j$.  Then there is a dart $x_i\in v_i$ such that $\{x_1,x_2,x_3\}$ is
a face of $H$.

\label{def:tame8}
\item There are never two nodes of type $(4,0,0)$ that are adjacent to one another.

\item
\label{def:tame3}
(Quadrilateral types) 
Let $v_1,v_2,v_3,v_4$ be any four distinct nodes of $H$ such that $v_i$ is adjacent
to $v_{i+1}$ for $i=1,2,3,4$ (setting $v_5=v_1$).  Then darts $x_i\in v_i$ can be chosen
so that the faces containing $x_i$ fall into one of the four patterns depicted
in Figure~\ref{fig:quadtype}.

\item 
\label{def:tame4}
The cardinality (degree) of a node of type $(p,q,r)$ is at most five if $r>0$.

\item 
\label{def:tame5}
The degree of each node of $H$ is at least $2$ and at most $6$.

\label{def:tame7}
\item There exists a Voronoi weight assignment of total weight at most $\mu(\Lambda_{dod})$.

\end{enumerate}
\end{definition}

\subsection{Classification}\label{sec:class}

Two hypermaps $(D,e,n,f)$ and $(D',e',n',f')$ are properly isomorphic if there
is a bijection between $D$ and $D'$ that is equivariant for the face, node, and
edge permutations.  Each hypermap $(D,e,n,f)$ has a mirror image:
  $$
  (D,f n, n^{-1},f^{-1}).
  $$
  An improper
isomorphism between two hypermaps is a proper isomorphism between one hypermap and the
mirror image of the other.
Two hypermaps are isomorphic if there is a proper or improper isomorphism between
the two hypermaps. If a hypermap is involutive or planar, then so is every isomorphic
hypermap.

There is an archive of tame Voronoi hypermaps
\cite{McLaughlin:2008:KeplerCode}. This archive contains $206$
hypermaps.

\begin{theorem}\label{thm:class}  
If $H$ is a tame Voronoi hypermap, then it is isomorphic to
some hypermap in the archive.
\end{theorem}

The proof of this theorem relies on the piece of 
computer code described in Section~\ref{sec:gg}.  The reader is referred to that section
for a description of the details of the theorem.

\subsection{Counterexamples are tame Voronoi hypermaps}

The following theorem proves that every potential counterexample to the Dodecahedral
conjecture gives a tame Voronoi hypermap.  In particular, the classification
of tame Voronoi hypermaps gives an explicit case enumeration of the possible combinatorial structures of a counterexample.  The following section will eliminate each case in the enumeration.
This will eliminate all possible counterexamples to the Dodecahedral conjecture.

\begin{theorem}  Assume that $\Lambda$ is a counterexample to the Dodecahedral
conjecture.  Without loss of generality,  assume that $0\in\Lambda=\Lambda(0,2t_{dod})$;
that the cardinality of $\Lambda^*$ is at least $13$; and that $G(\Lambda)$ is biconnected.  Let $H=(D,e,n,f)$ be the hypermap attached to the graph $G(\Lambda)$.
 For every dart $x\in D$ in face $F$, set
$$w(x) = \mu(\Lambda,U_F).$$
Then $H$ is a tame Voronoi hypermap and $w$ is a Voronoi weight assignment on $H$ of total weight at most $\mu(\Lambda_{dod})$.
\end{theorem}

Since the definition of tame Voronoi hypermap is a long enumeration of
different properties, the proof of this theorem breaks into a long
enumeration of lemmas, each establishing one property. The statement
of the theorem specifies the weight assignment $w$. The verification
that $w$ is a Voronoi weight, breaks into separate lemmas for each
property in the definition. The article \cite{Hales:2002:Dodec}
devotes many pages to the proofs of these lemmas. This articles
sketches the proofs and refers the reader to the fuller version for
details. Turn to an item by item discussion of the properties. The
first several are elementary.

\subsubsection{Involutive}
{\it The hypermap $H$ is involutive and planar.}  This has been established in Section~\ref{sec:ph}.

\subsubsection{Connected}
{\it The set of darts $D$  is a single orbit under $\tangle{e,n,f}$.}  This is also
contained in Section~\ref{sec:ph}. It follows directly from the
connectedness of $G(\Lambda)$.

\subsubsection{Simple}

{\it Every face of $H$ meets every node of $H$ in at most
dart.}  This is a direct consequence of the biconnectedness of $G(\Lambda)$.
Indeed, if by following the face permutation $x,f x, f^2 x,\ldots, f^r x= y$,
the darts $x\ne y$ lie at the same node $v$; then $v$ is an articulation vertex
of $G(\Lambda)$, and the graph is not biconnected.

\subsubsection{Node cardinality}

{\it The number of nodes is at least $13$.}  This is a consequence of the assumption
that the cardinality of $\Lambda^*$ is at least $13$ and the bijection in
Section~\ref{sec:ph} between nodes of the hypermap and $\Lambda^*$.

\subsubsection{Face cardinality}\label{sec:face}

{\it The cardinality of each face of $H$ is at least $3$ and at most $7$.}
By  definition, the face map on a dart $x$ takes the form
\begin{equation}\label{eqn:fx}
x = (0,v,w,u),\quad
f x = (0,w,\sigma_w^{-1}v,v),\quad
f^{-1} x = (0,u,v,\sigma_u v).
\end{equation}
In a biconnected graph with more than two vertices,
every vertex has degree at least two.
Thus every node of the hypermap has degree
at least two.  Thus, $\sigma_u v\ne v$ and $f x\ne f^{-1} x$.
Also, $w,u\in E(v)$, which does not contain $v$.  Thus, 
the form of $x,fx$, and 
$f^{-1} x$
in (\ref{eqn:fx}) shows these darts are distinct, and
the
face contains at least three distinct darts.

If some face of $H$ has cardinality at least $8$, then by Theorem~\ref{thm:main}
$$
\mu(\Lambda) =\sum_{U\in [Y(\Lambda)]}\mu(\Lambda,U) > t_8 =\mu(\Lambda_{dod}).
$$
Thus, $\Lambda$ is not a counterexample, as was assumed.

\subsubsection{Triangle types}

{\it Let $v_1,v_2,v_3$ be any three nodes of $H$ such that $v_i$ is
adjacent to $v_j$ for each $i\ne j$. Then there is a dart $x_i\in v_i$
such that $\{x_1,x_2,x_3\}$ is a face of $H$.}
This is a restatement of Lemma~\ref{lemma:enclosed} in terms of the
combinatorial properties of hypermaps.

\subsubsection{Adjacent degrees}

\begin{figure}[htb]
  \begin{center}
    \includegraphics[scale=0.50]{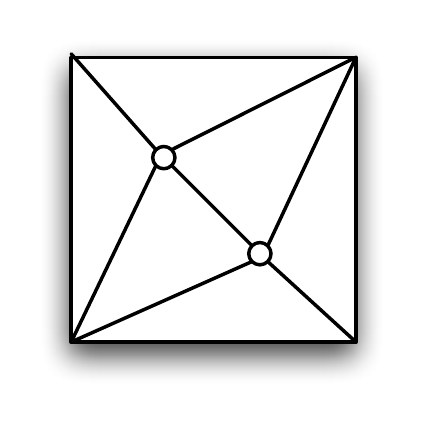}
   \end{center}
  \caption{Two adjacent vertices of type $(4,0,0)$}
  \label{fig:adj4}
\end{figure}

{\it There are never two nodes of type $(4,0,0)$ that are adjacent to one another.}
If there are two adjacent nodes of type $(4,0,0)$, then the graph takes the
shape of Figure~\ref{fig:adj4}.  This is an impossible configuration in a packing $\Lambda$
for purely geometric reasons.  It has nothing to do with the value of $\mu(\Lambda)$
and volumes of truncated Voronoi cells.
The impossibility proof appears as \cite[Lemma~3.8]{Hales:2002:Dodec}. It is a reparametrization
of \cite[Prop.4.2]{Hales:1997:DCG}.  This is one of the most delicate reparametrizations.

\subsubsection{Quadrilateral types}

{\it Let $v_1,v_2,v_3,v_4$ be any four distinct nodes of $H$ such that $v_i$ is adjacent
to $v_{i+1}$ for $i=1,2,3,4$ (setting $v_5=v_1$).  Then darts $x_i\in v_i$ can be chosen
so that the faces containing $x_i$ fall into one of the four patterns depicted
in Figure~\ref{fig:quadtype}.}
Let $G'$ be the spherical graph on the vertex set $\{v_1,v_2,v_3,v_4\}$
with edges $\{v_i,v_{i+1}\}$ for $i=1,2,3,4$.  Then by
Jordan curve theorem for polygons, $Y(G')$ consists of two connected components.
The result of \cite[Lemma~3.8]{Hales:2002:Dodec} cited in the previous proof states more precisely
that exactly one connected component $U$ of $Y(G')$ has solid angle less than $2\pi$,
and that $U\cap\Lambda$ contains at most one point.
If $U\cap\Lambda$ is empty, then the first pattern of Figure~\ref{fig:quadtype} occurs.

\begin{figure}[htb]
  \begin{center}
  \includegraphics[scale=0.50]{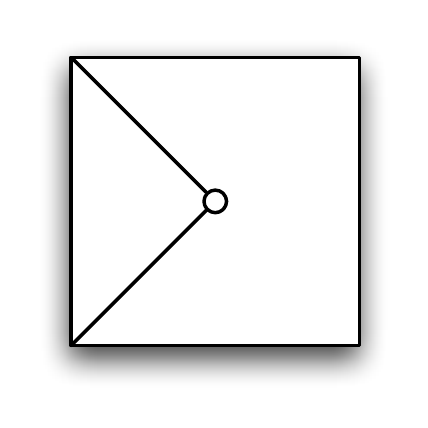}
  \end{center}
  \caption{A quad type that is excluded}
  \label{fig:tripent}
\end{figure}

In remaining cases,  $U\cap\Lambda$ contains a single point $v_0\in\Lambda$.  Again, by
the Jordan curve theorem and the planarity of $G(\Lambda)$, 
all edges $\{v_0,u\}$ in $G(\Lambda)$ have the form
$\{v_0,v_i\}$ for $i=1,2,3,4$.  Section~\ref{sec:face} shows that each node has  degree
at least $2$.  Thus,  $v_0$ is adjacent to $2$, $3$, or $4$ of the vertices $v_i$.
Figure~\ref{fig:quadtype} gives all such connection patterns of $v_0$ with $v_i$, except
the one shown in Figure~\ref{fig:tripent}.  

Thus, it is enough to show that pattern of Figure~\ref{fig:tripent} does not occur 
in any counterexample to the Dodecahedral conjecture.  This graph contains a triangle
$F$ and a pentagon $F'$, with corresponding standard components
$U=U_F$ and $U'=U_{F'}\in[Y(\Lambda)]$.
An important estimate \cite[Lemma~10.1]{Hales:2002:Dodec} gives that
$$
\mu(\Lambda,U) + \mu(\Lambda,U') > 0.168.
$$
If there is some other face  with $n\ge 4$ sides, then 
Theorem~\ref{thm:main} gives
$$
\mu(\Lambda) \ge \mu(\Lambda,U)+\mu(\Lambda,U') + t_n \ge 0.168 + 0.031 > \mu(\Lambda_{dod}).
$$
Otherwise, pick any four vertices in $\Lambda^*$ other than $v_0,v_1,\ldots,v_4$.  Let $U_1,\ldots,U_r$ be the standard components
indexed by triangles of $G(\Lambda)$
that contain one of these four vertices.  
Another estimate \cite[Lemma~5.2]{Hales:2002:Dodec} gives
$$
\sum_{i=1}^r \mu(\Lambda,U_i) \ge 4 (0.004) \text{ and hence }
\mu(\Lambda) \ge 0.168 + 4(0.004) > \mu(\Lambda_{dod}).
$$
This shows that $\Lambda$ is not a counterexample.

\subsubsection{Degree}

{\it The degree of a node of type $(p,q,r)$ is at most five if $r>0$.}
Indeed, it is a direct consequence of the inequalities of \cite{McLaughlin:2008:KeplerCode} that if $U_1,\ldots,U_n$, for $n=p+q+r\ge 6$, are the standard components of the faces around the node, then 
  $$
  \sum_{i=1}^n \mu(\Lambda,U_i) > \mu(\Lambda_{dod}).
  $$


\subsubsection{Degrees}\label{sec:degrees}

{\it The degree of each node of $H$ is at least $2$ and at most $6$.}
Section~\ref{sec:face} has already shown that the degrees are at least $2$.
Let $(p,q,r)$ be the type of node $v$.  If $r>0$, then the previous
property bounds the degree at five.  Assume $r=0$.  The
proof in this case is deferred until Section~\ref{sec:pq}.

\subsubsection{Total weight}

{\it There exists a Voronoi weight assignment of total weight at most $\mu(\Lambda_{dod})$.}
By definition, a counterexample $\Lambda$ is a packing such that
\begin{equation}\label{eqn:178}
\mu(\Lambda)=\sum_{U\in [Y(\Lambda)]} \mu(\Lambda,U) \le \mu(\Lambda_{dod}).
\end{equation}
The set  $[Y(\Lambda)]$ of standard components is in bijection with the faces of $H$.
By the definition of $w$ in the statement of the theorem, (\ref{eqn:178}) can
be rewritten as
$$
\sum_{x\in[D/f]} w(x) \le \mu(\Lambda_{dod}).
$$
This is exactly what it means for $w$ to have total weight at most $\mu(\Lambda_{dod})$.

\subsubsection{Voronoi weight}\label{sec:pq}

{\it The weight assignment $w$ is a Voronoi weight.}
This can be expressed as the following claims.
\begin{enumerate}
\item If the face $F$  has cardinality $n\ge3$ then $\mu(\Lambda,U_F)> t_n$.
\item Let $F\subset D$ 
be any face with cardinality $m \ge 5$.  Let $y\in F$.
Let $V$ be a set of vertices of $F$, such
that no two are adjacent to one another.  
Assume that the type of each node of $V$ is $(4,0,1)$.
Let $X$ be the set of triangles at the vertices $V$.
that is, the set of darts in nodes in $V$ except those in $F$.
Let $m'$ be the cardinality of $V$.
Then
\begin{equation}\label{eqn:m'}
\mu(\Lambda,U_F) + \sum_{F'\in X} \mu(\Lambda,U_{F'}) > t_m  +  0.016\, m'.
\end{equation}
\item If the node of $x$ has type $(p,q,0)$ and degree $m=p+q$, then
  \begin{equation}\label{eqn:pq}
  \sum_{i=1}^m \mu(\Lambda,U_i) > b(p,q).
  \end{equation}
\end{enumerate}

The first of these claims is a direct consequence of
Theorem~\ref{thm:main}. Consider the second claim. Form the graph
$G'(\Lambda)$ and partition $U_F$ into subcomponents by the algorithm
described in Section~\ref{sec:sub}. A result in Tarski arithmetic
states that two internal pairs $\{u,w\}$ and $\{u',w'\}$ with
$|u-w|,|u'-w'|\le\sqrt8$ do not cross each other in the sense of
Section~\ref{sec:sub}~\cite[Sec.~7.2.1]{Hales:2002:Dodec}. It
follows that every internal pair $\{u,w\}$ such that $|u-w|\le\sqrt8$
forms an edge of the graph $G'$. The bound on $t_m$ is obtained by
superadditivity, from a collection of inequalities
$\mu(\Lambda,U_i)>D_{dod}(n,k)$ for each subcomponent.

Let $v$  be a node in $V$.  Let $(0,v,u,w)$ be the dart in $F$
at node $v$.  Let $U_v$ be the subcomponent indexed by the
triangle $\{v,u,w\}$ in $G'$.
At a node of type $(4,0,1)$, by angle considerations,
the pair $\{u,w\}$ is necessarily internal.
Let $U'_1,\ldots,U'_4$ be the standard components indexed
by the four triangles at $v$. 
By summing over $V$, the following two inequalities
imply (\ref{eqn:m'}).
\begin{enumerate}
\item If $|u-w|>\sqrt8$, then
  $$
  \sum_{i=1}^4\mu(\Lambda,U'_i) > 0.016.
  $$
\item If $|u-w|\le\sqrt8$, then
  $$
  \mu(\Lambda,U_v)+ \sum_{i=1}^4\mu(\Lambda,U'_i) > D_{dod}(3,1) + 0.016.
  $$
\end{enumerate}
These inequalities are obtained by summing over interval arithmetic
inequalities for each term $\mu(\Lambda,\cdot)$.   
A detailed proof appears at \cite[Theorem~8.1]{Hales:2002:Dodec}.

Consider the final claim.  
The argument at \cite[p.14]{Hales:2002:Dodec} goes as follows.
Let $v\in\Lambda^*$ be a vertex.
To have type $(p,q,0)$ means that there are $p$ triangles and $q$ quadrilaterals
in the graph $G(\Lambda)$ at the vertex $v$.
Interval arithmetic can be used to compute lower and upper bounds on the
azimuth angles $\op{azim}(U_F,v)$ when $U_F$ is a triangle or quadrilateral.
These bounds are \cite[F.2.1,F.4]{Hales:2002:Dodec}
$$
\begin{array}{lrlll}
\text{triangle:} &  0.856147 &< \op{azim}(U,v) &< 1.88673,\\
\text{quadrilateral:} & 1.15242 &< \op{azim}(U,v) &< 3.25887.\\
\end{array}
$$
\begin{equation}\label{eqn:2pi}
\sum^m \op{azim}(U_i,v) = 2\pi.
\end{equation}
By the bounds on the azimuth angles, this equality can only be satisfied for
special $(p,q)$. Anything else is a geometric impossibility.  The feasible pairs $(p,q)$
are listed in \cite[Lemma~6.1]{Hales:2002:Dodec}.  According to this list, $p\le 7$ and $q\le 5$.

Interval arithmetic methods establish a list of nonlinear inequalities
relating $\op{azim}(U_F,v)$ to $\mu(\Lambda,U_F)$ when $F$ is a
triangle or quadrilateral \cite{McLaughlin:2008:KeplerCode}. If free variables
$\optt{azim}(F)$ and $\optt{mu(F)}$ are substituted into these
inequalities for $\op{azim}(U_F,v)$ and $\mu(\Lambda,U_F)$, then the
resulting inequalities are linear in these variables. A system of
linear inequalities results. (This is the method of linear
relaxation.) A lower bound on the left-hand side of (\ref{eqn:pq}) is
the solution of the linear program
$$
\min \sum_F \optt{mu(F)}
$$
subject to 
$$
\sum_F \optt{azim(F)} = 2\pi,
$$
and to the system of linear inequalities.  The linear program is run for
each $p,q$ and a constant $b(p,q)$ slightly smaller than the minimization was picked.
This gives the table of values.

If $\Lambda$ is a counterexample to the Dodecahedral conjecture with a vertex
of type $(p,q,0)$, then
$$\mu(\Lambda_{dod}) \ge \mu(\Lambda) > b(p,q).$$
An inspection of the list of constants $b(p,q)$ in the table shows that that this implies
that $p+q\le 6$.  That is, the degree of every vertex of type $(p,q,0)$
is at most $6$.  This
is the property needed in Section~\ref{sec:degrees}.

\section{Linear Programs}
\label{sec:linear-programs}

This section discusses Theorem~\ref{thm:graph-system}.  It is one of
the main steps in the proof of the Dodecahedral conjecture.
The discussion begins with the terminology used in the
statement of the theorem.

\begin{definition} A hypermap system is a pair $(H,\Phi)$,
where $H=(D,e,n,f)$ is a hypermap, and $\Phi$ is a finite set constraints on $H$.  More precisely, let $V$ be the vector space of
real-valued functions on $D$.  Each $\phi\in\Phi$ is a boolean valued function 
$\phi:V^\ell\to \{\op{true},\op{false}\}$
for some $\ell$. (Assume $\ell$ is independent of $\phi\in \Phi$).

The hypermap system $(H,\Phi)$ is said to be feasible, if
there is some $x=(x_1,\ldots,x_\ell)\in V^\ell$ such that
$\phi(x)$ holds for all $\phi\in\Phi$. Otherwise,  
the system is infeasible.
\end{definition}

In our hypermap systems, the constraints $\Phi$ 
are generated from a list of
about one hundred generic parametrized constraints, with a parameter
running over the darts in a hypermap. 
For example, a parameterized constraint $c(H,x)$ might be 
interpreted as an inequality
relating  $\mu(\Lambda,U_F)$, 
the solid angle of $U_F$, and
azimuth angles, whenever $F$ is a triangular  face containing the dart $x$.
As $x$ runs over the darts in a hypermap $H$, 
the generic constraint $c(H,x)$ 
evaluates to a constraint for each dart $x$ in every triangular face of $H$, 
yielding dozens of particular constraints for $H$.
  We do not list the generic constraints here, but they appear
in~\cite{Hales:2002:Dodec} and also
at~\cite{McLaughlin:2008:KeplerCode}.

For each tame Voronoi hypermap $H$, the generic constraints specialize
to a particular set of constraints $\Phi_{dod}(H)$. Each set $\Phi_{dod}(H)$
contains about $5000$ constraints. The hypermap
systems $(H, \Phi_{dod}(H))$ arising in this way are called 
\emph{Voronoi hypermap systems}.

\begin{theorem}\label{thm:graph-system}  Let 
$H$ be any tame Voronoi hypermap. 
Then $(H,\Phi_{dod}(H))$ is infeasible.
\end{theorem}

This proof is carried out by computer as a collection of linear
programs. This is one of the three major parts of the proof of the
Dodecahedral conjecture that have been carried out by computer. This
article describes the relationship between the feasibility of
$(H,\Phi)$ and a linear programming feasibility problem. It also
describes some details of the implementation of the code.

Section~\ref{sec:class} enumerates of all tame Voronoi
hypermaps.  Thus, the proof of
Theorem~\ref{thm:graph-system} may proceed case by case.
Let's focus attention for a moment on one tame Voronoi hypermap
$H=(D,e,n,f)$ and the corresponding Voronoi hypermap system
$(H,\Phi)$. A simple strategy will show that it is infeasible. For
some $\ell\in\ring{N}$, each constraint $\phi$ is a function on
$V^\ell$, where $V$ is the vector space of real-valued functions on
$D$. Thus, $V^\ell$ can be identified with $\ring{R}^m$, where $m=
\ell\, \#(D)$. An inspection of the form of the generic constraints
$\phi\in \Phi$ described above reveals they all have a very special
form. They are all linear constraints on $\ring{R}^m$.

Some of the linear constraints
carry guard conditions.  That is, some constraints have the form
  \begin{equation}\label{eqn:guard}
  (A x < b)  \Rightarrow (A' x \le b'),
  \end{equation}
for $x\in\ring{R}^m$, and various matrices $A,A'$ and vectors
$b,b'$.  (The vector inequality $a \le b$ means
that $a_i\le b_i$ for every component of the vectors $a,b$.)
The constraint $(A x < b)$ is called a guard condition.
Variations are allowed in which some of the inequalities in the
guard condition are weak and some of the inequalities in the
consequent are strict.

The collection of all inequalities that do not have a guard 
condition is a system of linear inequalities.  Standard linear
programming packages can be used to determine whether this
system of linear inequalities has a feasible solution.  If this
linear program is infeasible, then the hypermap system $(H,\Phi)$
is clearly also infeasible.  When this happens, a
proof of the infeasibility of $(H,\Phi)$ results.

When this fails, the constraints with guard conditions are
used.
The introduction of a constraint that has a nontrivial guard condition
involves multiple steps.  
The constraint (\ref{eqn:guard}) can be rewritten in logically
equivalent form as
  $$
   (A_{1} x \ge b_{1}) \lor \cdots \lor
   (A_{r} x \ge b_{r}) \lor (A' x \le b'),
  $$
where $A_{i}$ and $b_{i}$ are the rows of $A$ and $b$.
Taking each disjunct in turn, one linear inequality at a time
is added to the system
of linear inequalities, and the resulting system is shown to be
infeasible.  When each
systems are infeasible, then $(H,\Phi)$ itself is infeasible.

This discussion may give the impression 
that a great many linear programming feasibility
problems are created in this manner.  In practice, nearly all
of the hypermap systems are eliminated in the first pass, without
requiring recourse to the guard conditions.  

\subsection{Counterexample implies feasibility}

\begin{theorem}\label{thm:feasible}  Let $\Lambda$ be a counterexample to the Dodecahedral conjecture.  Let $(H,\Phi)$ be the Voronoi hypermap system attached to $\Lambda$.
Then $(H,\Phi)$ is feasible.
\end{theorem}

\begin{proof} It is enough to give some assignment of
$x=(x_1,\ldots,x_\ell)\in V^\ell$ 
such that $\phi(x)$ holds for all $\phi\in\Phi$.
Each $x_i$ is a real-valued function on darts.
The notation  for the functions $x_i$ 
has been set up in a way that suggests the assignment.  
Let $\alpha = (0,v,u,w)$ be a dart.  Let $F$ be the face containing $\alpha$. Make the following settings.  The first row of this table gives the functions $x_i$ in the notation of~\cite{McLaughlin:2008:KeplerCode}. The second row gives the real number $x_i(\alpha)$.
$$
\begin{array}{llllllllllll}
\optt{yn} & \optt{ye} & \optt{sol} & \optt{azim} & \optt{mu} & \optt{omega} \\
|v|               & |v-u|           & \sol(U_F)   & \op{azim}(U_F,v) & \mu(\Lambda,U_F) & \omega(\Lambda,U_F).\\
\end{array}
$$

All of the predicates $\Phi$ can be shown to hold for this assignment,
either as a consequence of definitions, as consequences of geometrical
facts, or as consequences of interval arithmetic calculations. 

The interval arithmetic calculations for the predicates with guard
conditions were verified in two stages. The polygonal face is
triangulated, and the standard component $U_F$ is partitioned in a
corresponding way into parts $U_1,\ldots,U_k$. In the first stage, the
guard hypotheses were used to obtain interval arithmetic bounds on the
lengths of the internal edges of the triangulation. In the second
stage, these edge length bounds are used as hypotheses in further
interval arithmetic bounds that give lower bounds on
$\omega(\Lambda,U_i) > c_i$. Then
  $$\omega(\Lambda,U_F) = \sum_{i=1}^k\omega(\Lambda,U_i) >\sum_{i=1}^k c_i=c.$$
This becomes a general predicate on darts of the form:

$$
\optt{omega}(\alpha) > c.
$$
\end{proof}

Since no tame Voronoi hypermap is feasible, by the argument of
Section~\ref{sec:proof-outline} the Dodecahedral conjecture is
established.

\subsection*{Thanks}

This research was supported by NSF grant 0503447. Many of the
nonlinear inequalities for the proof of the Dodecahedral conjecture
were proved by computer by S. Ferguson. We wish to thank him for his
generous support of this project. We wish to thank C. Anghel for
helpful comments.


\bibliographystyle{abbrv}
\bibliography{dodec}


\end{document}

%% file: macros.tex
\def\displayallproof{f} 

\def\swallowed{\relax}
\def\swallow#1\swallowed{}
\newenvironment{iproved}{}{}

\def\hideproof{\renewenvironment{iproved}{%
   \centerline{\it -- Proof Proofed --}
  \renewenvironment{itemize}{}{}
  \renewenvironment{enumerate}{}{}
  \def\item{\relax}
  \catcode13=12
  \swallow
}{}}
\def\showproof{\renewenvironment{iproved}{\begin{proof}}{\end{proof}}}
\def\resetproved{\if\displayallproof t\showproof\else\hideproof\fi}

\def\rating#1{\if\verbose t{{\textsc {[rating={\ensuremath {#1}}].\ }}}\else{}\fi}
\def\formal#1{\if\verbose t{{\tt [formal:#1].\ }}\else{}\fi}
\def\footformal#1{\if\verbose t{\footnote{\formal{#1}}}\else{}\fi}
\def\usage#1{\if\verbose t{\symbolfootnote[1]{#1}}\else{}\fi}

\def\FIXX#1{\if\verbose t{\footnote{\tt [#1]}}\else{}\fi}

\def\dcg#1#2{{\if\verbose t{{\tt{[DCG-#1]}}\index{References}{ZC{#2 #1}@{DCG-#1}|page{#2}}}\else{}\fi}}
\def\tlabel#1{\label{#1}\if\verbose t{{\tt [#1].\ }%
   \index{References}{#1|itt}}\else{}\fi}
\def\ifverbose#1{\if\verbose t{{#1}}\else{}\fi}
\def\skipXX#1{}

\def\guid#1{{\tt[#1].\ }\index{References}{ZA{#1}@{#1}|itt}}

\def\tarfe#1#2{{\foote{{\tt[#1]}~\csname #1-sum\endcsname%
  \expandafter\ifx \csname #1-guid\endcsname \relax [XX-BAD-IDENTIFIER]\else{}\fi
  }
  {\csname #1-guid\endcsname{#2}}}}

\def\tarfE#1{\tarfe{EE}{-\csname #1\endcsname}}
\newtoks\mysummary
\newtoks\myguid
\newtoks\myname
\def\separator{\if\tarskipagesep t{\clearpage}\else%
  \smallskip\hrule height 0.5pt depth 0pt width 60pt\bigskip\fi}
\def\identity#1{#1}
\def\foote#1#2{\footnote{#1~{\it (#2)}}}
\def\marker#1#2{\par\hangafter1\hangindent=1em\noindent{\bf #1:}\quad #2}

\def\summary#1{\marker{Summary}#1\mysummary={#1}
  \if\tarskistrip t{
  \def\sqr{sqrt}
  \def\sqrt{sqrt}
  \def\op{\identity}
  \def\epsilon{epsilon}
  \def\Delta{Delta}
  \def\ups{upsilon}
  \def\beta{beta}
  \def\rho{rho}
  \def\CalE{E}
  \def\chi{chi}
  \def\eta{eta}
  \write\myhtml{  "#1", // summary}
  }
  \else\relax\fi
}
\def\guid#1{\marker{ID}#1
   \index{guid}{#1}
   \immediate\write\mywrite{\global\expandafter\def\csname \the\myname-sum \endcsname{\the\mysummary}}
   \immediate\write\mywrite{\global\expandafter\def\csname \the\myname-guid \endcsname{#1}}
   \immediate\write\myhtml{  "#1" // guid}
}

\long\def\symbolfootnote[#1]#2{\begingroup%
\def\thefootnote{\ensuremath{\fnsymbol{footnote}}}\footnote[#1]{#2}\endgroup}
\def\hide#1{}

 \DeclareMathOperator{\sol}{sol}
 
\def\CalE{{\mathcal E}}

\def\CalS{{\mathcal S}}

\newcommand{\ring}[1]{\mathbb{#1}}

\def\op#1{{\operatorname{#1}}}
\def\optt#1{{\operatorname{{\texttt{#1}}}}}

\def\azim{\op{azim}}
\def\sol{\operatorname{sol}}

\def\quo{\operatorname{quo}}

\def\tangle#1{\langle #1\rangle}

\def\ups{\upsilonup} 
\def\mid{\ :\ }

\def\sqr{\sqrt}

\def\=#1{\accent"16 #1}

\expandafter\def\csname 2t0-doesnt-pass-through\endcsname{01}
\expandafter\def\csname E1453\endcsname{02}   
\expandafter\def\csname dcg-p89\endcsname{03}   
\expandafter\def\csname enclosed-v\endcsname{04}   
\expandafter\def\csname 277\endcsname{05}   
\expandafter\def\csname 245\endcsname{06}   
\expandafter\def\csname 245bis\endcsname{07}   
\expandafter\def\csname node\endcsname{08}   
\expandafter\def\csname E:part4:1\endcsname{09}   
\expandafter\def\csname E:part4:2\endcsname{10}   
\expandafter\def\csname E:part4:3\endcsname{11}   
\expandafter\def\csname E:part4:4\endcsname{12}   
\expandafter\def\csname E:part4:5\endcsname{13}   
\expandafter\def\csname E:part4:6\endcsname{14}   
\expandafter\def\csname E:part4:7\endcsname{15}   
\expandafter\def\csname E:part4:8\endcsname{16}   
\expandafter\def\csname E:part4:9\endcsname{17}   
\expandafter\def\csname E:part4:10\endcsname{18}   
\expandafter\def\csname dcg-p142\endcsname{19}   
\expandafter\def\csname convex-quad\endcsname{20}   
\expandafter\def\csname last:E\endcsname{21}
\expandafter\def\csname rem2.7\endcsname{22}     
\expandafter\def\csname pass-anchor\endcsname{23}